\documentclass[11pt]{amsart}
\usepackage{fullpage}
\usepackage[utf8]{inputenc}
\usepackage[T1]{fontenc}
\usepackage{graphicx}
\usepackage{amsmath}
\usepackage{amsfonts}
\usepackage{graphicx}
\usepackage{dsfont}
\usepackage{color}

\usepackage{subfigure}

\usepackage[pdftex]{hyperref}
\usepackage{pdfpages}

\newtheorem{theo}{Theorem}[section]
\newtheorem{cor}[theo]{Corollary}
\newtheorem{rem}[theo]{Remark}

\newtheorem{propo}[theo]{Proposition}
\newtheorem{lemme}[theo]{Lemma}

\newcommand{\E}{\mathbb{E}}
\newcommand{\R}{\mathbb{R}}

\newcommand{\F}{\mathcal{F}}

\newcommand{\N}{\mathbb{N}}
\newcommand{\dd}{\mathrm{d}}
\newcommand{\II}{\mathcal{I}}
\newcommand{\PPP}{\mathcal{P}}

\numberwithin{equation}{section}

\title[Monte-Carlo semi-lagrangian scheme]{Analysis of the Monte-Carlo error in a hybrid semi-lagrangian scheme}
\author[C.-E. Br\'ehier]{Charles-Edouard Br\'ehier}
\address{INRIA \& ENS Cachan Bretagne - IRMAR, Universit\'e Rennes 1\\Avenue Robert Schuman\\F-35170 Bruz\\France}
\email{charles-edouard.brehier@bretagne.ens-cachan.fr}
\author[E. Faou]{Erwan Faou}
\address{INRIA \& ENS Cachan Bretagne - IRMAR, Universit\'e Rennes 1\\Avenue Robert Schuman\\F-35170 Bruz\\France}
\email{Erwan.Faou@inria.fr}

\begin{document}

\begin{abstract}
We consider Monte-Carlo discretizations of partial differential equations based on a combination of semi-lagrangian schemes and probabilistic representations of the solutions. 
We study  the Monte-Carlo error in a simple case, and show that under an anti-CFL condition on the time-step $\delta t$ and on the mesh size $\delta x$ and for $N$ - the number of realizations -  reasonably  large, we control this error by a term of order $\mathcal{O}(\sqrt{\delta t /N})$. 
We also provide some numerical experiments to confirm the error estimate, and to expose some examples of equations which can be treated by the numerical method.
\end{abstract}

\maketitle

%\part[Analysis of the Monte-Carlo error in a hybrid semi-lagrangian scheme]{Analysis of the Monte-Carlo error in a hybrid semi-lagrangian scheme}

%\vspace{3cm}

\section{Introduction}

The goal of this paper is to analyze and give some error estimates for numerical schemes combining the principles of semi-lagrangian and Monte-Carlo methods for some partial differential equations. 
Let us describe the method in a very general case by considering first a linear transport equation of the form 
$$
\partial_t u(t,x) = f(x) \cdot \nabla u(t,x), \quad x \in  \R^d, \quad u(0,x) = u_0(x), 
$$
where $u_0$ is a given function. Under some regularity assumptions and existence of the flow associated with the vector field $f(x)$ in $\R^d$, the solution of this equation is given by the characteristics representation $u(t,x) = u(0, \varphi_{t}(x))$, where $\varphi_t(x)$ is the flow associated with the ordinary differential equation $\dot y = f(y)$ in $\R^d$. 
In this context, semi-Lagrangian schemes can be described as follows. Let us consider a set of grid nodes $x_j$, $j \in K$ in $\R^d$ ($K = \N$ or a finite set) and an interpolant operator  $\II$ mapping  vectors of values at the nodes, $(u_j) \in \R^K$ to a function $(\II u)(x)$ defined over the whole domain. In this paper we will consider the case where $x_j = j (\delta x)$, $j \in \mathbb{Z}^d$, $\delta x$ is the space mesh size, and $\II$ a standard linear interpolation operator. Given approximations $u_j^n$ of the exact solution $u(t_n,x_j)$ at times $t_n = n (\delta t)$ and points $x_j$, the previous formula gives an approximation scheme for $u_j^{n+1}$ obtained by solving the ordinary differential equation $\dot y = f(y)$ between $t_n$ and $t_{n+1}$: $u_j^{n+1} = (\II u^n)(\Phi_{\delta t}(x_j))$, where $\Phi_h$ is the numerical flow associated with a time integrator. 

These methods are particularly interesting when the vector field $f(x;u)$ depends on the solution $u$ making the transport equation nonlinear, see for instance \cite{Staniforth, Falcone, Sonnen} and the references therein. This is the case when an advection term is present for instance, or for Vlasov equations (see for instance \cite{crouze}). In these situations, standard semi-lagrangian schemes are based on solving equations of the form 
$$
\partial_t  u(t,x) = f(x;u^n) \cdot \nabla u (t,x), 
$$
between $t_n$ and $t_{n+1}$, where $u^n$ denotes the solution at time $t_n$. In other words the vector field is frozen in $u^n$ (in the language of geometric numerical integration, it is Crouch and Grossman method, see \cite{crouch}). If moreover the vector field $f(x;u)$ possesses some geometric structure for all functions $u$, the numerical integrator can be chosen to preserve this structure (for example symplectic integrator in the Vlasov case). %\red structure de la derni\`ere phrase: "geometric structure", "these?" \black

In many situations, a diffusion term is present, and the equation can be written (in the linear case)
\begin{equation}
\label{defPDE0}
\partial_t u(t,x) = \textstyle\frac{1}{2} \sigma(x) \sigma(x)^T \Delta u + f(x) \cdot \nabla u(t,x), \quad x \in \R^d, \quad u(0,x) = u_0(x), 
\end{equation}
where $\sigma(x)$ is a $d \times n$ matrix. 
In this case the solution admits the characteristic representation 
$$
u(t,x) = \E \, u_0(X_t^x),
$$
where $X_t^x$ is the stochastic process associated with the stochastic differential equation 
\begin{equation}
\label{defEDS0}
\dd X_{t}^{x}=f(X_{t}^{x})\dd t+ \sigma(X_t^x) \dd B_t,\quad X_{0}^{x}=x,
\end{equation}
where $(B_{t})_{t\geq 0}$ is a standard $n$-dimensional Brownian Motion.
%\black 
%
%We consider a class of numerical methods which combine the principles of semi-lagrangian and Monte-Carlo methods. We present their construction and their analysis in a simple case: we consider a diffusion equation on the domain $D=\mathbb{R}^d$
%\begin{equation}\label{defEDP}
%\begin{gathered}
%\frac{\partial u(t,x)}{\partial t}=\frac{1}{2}\Delta u(t,x)+c(x).\nabla u(t,x), \text{ for any } 0<t\leq T \text{ and }x\in\R^d\\
%u(0,x)=u_0(x), \text{ for any } x\in\R^d,
%\end{gathered}
%\end{equation}
%with a given initial condition $u_0$, assumed to be smooth. The function $c$ is assumed to be smooth. We notice that with $D=\mathbb{R}^d$ there is no boundary condition, so that the presentation is simplified. The theoretical part of this article deals with the case with $c=0$, and periodic boundary conditions on the domain $D=(0,1)$. % In Chapter $5$, we give elements for more general situations.
%%with periodic boundary conditions on the domain $D=(0,1)$, and a smooth, periodic initial condition $u_0:\mathbb{R}\rightarrow\mathbb{R}$.
%
%
%
%The infinitesimal generator of this diffusion process is $\mathcal{L}=c.\nabla+\frac{1}{2}\Delta $, and we see that $\frac{\partial u}{\partial t}=\mathcal{L}u$.

In general, the law of the random variable $X_{t}^{x}$ is not explicitly known, and we are not able to compute the expectation. The classical approximation procedures for such problems are Monte-Carlo methods: 
if we assume that we are able to compute $N$ independent realizations $(X_{t}^{x,m})_{1\leq m\leq N}$ of the law of $X_{t}^{x}$, we can approach $u(t,x)$ with
\begin{equation}
\label{MCnaif}
\frac{1}{N}\sum_{m=1}^{N}u_0(X_{t}^{x,m}).
\end{equation}
In general, the variance of the random variables $u_0(X_{t}^{x,m})$ is of size $t$ and the law of large numbers ensures that the statistic error made is typically of order $\mathcal{O}(\sqrt{T/N})$ for an integration over the interval  $[0,T]$. To this error must be added the error in the approximation of the process $X_{t}^{x}$ by numerical schemes of Euler type for instance. This is error is of order $\mathcal{O}(\delta t )$, see for instance \cite{kloeden-platten, Mil-Tre, Talay96} and the reference therein for analysis of the weak-error in the numerical approximation of stochastic differential equations. 
If a global knowledge of the solution is required, the above operation must be repeated for different values of $x_j$
 on the grid. 
%
%
%We see that these discretization methods can be used, but that the associated variances in the Monte-Carlo method are not small, even if $\delta t$ and $\delta x$ are small. Moreover, for each new computed value the procedure must be started again. The method we now propose uses more carefully the time and space dependence of the solutions, and aims at reducing the variances. Moreover, we directly obtain a global approximation defined on the whole spatial domain, in an intrinsic way.
 
The numerical method we study in this paper is based on the Markov property of the associated stochastic processes: we have for any $x_j$ on the spacial grid and locally in time
\begin{equation}
\label{reprMarkov}
u(t_{n+1},x_j)=\E u(t_n,X_{\delta t}^{x_j}),%=\E u(t_1,x+\sqrt{t_2-t_1}\mathcal{N}),
\end{equation}
%where $\mathcal{N}$ is a standard Gaussian random variable.
which is the formula we aim at discretizing. 
Using the Euler method to compute a numerical approximation of $X_{\delta t}^{x_j}$, we end up with the following numerical scheme
\begin{equation}
\label{eq:schema}
u_{j}^{n+1}:=\frac{1}{N}\sum_{m=1}^{N}\bigl(\mathcal{I}u^n\bigr)(x_j + \delta t f(x_j)+\sqrt{\delta t}\mathcal{N}^{n,m,j}),
\end{equation}
where the random variables $(\mathcal{N}^{n,m,j})_{1\leq m\leq M}$ are independent standard Gaussian random variables. 
Note that the main difference between the standard Mont-Carlo method is that the average is computed at every time step. %Moreover a fundamental assumption is that independence with respect to the parameters $n$ and $j$ holds. On the one hand, we then break a symmetry of the problem in the space variable with a noise which is white in space; on the other hand, the noise has an effect on the convergence which is essential in the proof of our main Theorem \ref{TheTheorem}.

%We can remark that $v^n=\E u^n$.

Such numerical method were already introduced in \cite{faou}. 
The principle of using (random) characteristic curves in \eqref{reprMarkov} over a time interval of size $\delta t$ and to use an interpolation procedure to get functions defined on the whole domain fits in the semi-lagrangian framework. The addition of the Monte-Carlo approximation then justifies the use of the hybrid terminology.

As in the deterministic case described above, it is clear the method can be adapted in situations where the drift term $f(x)$ and the noise term $\sigma(x)$ depend on the solution $u$.  We will present in the end of the paper some numerical experiment in such nonlinear situations. 

Another remark is that different kind of boundary conditions can be considered. If the presentation made above was concerning the situation where the equation is set on $\R^d$, representation formulae such as \eqref{reprMarkov} hold true in the case of periodic boundary conditions, Dirichlet or Neumann condition on bounded domains. Again, we give some numerical examples in the end of the paper. 

But the main aim of this paper is to perform the numerical analysis of the scheme \eqref{eq:schema} in the simplest situation, that is $f = 0$, in dimension $d = 1$, with $\sigma(x) = 1$ and periodic intial conditions such that the transport equation is set on a domain  $D=(0,1) \subset \R$ with periodic boundary conditions. Note that using splitting strategy, this is the term that is new in \eqref{eq:schema} in comparison with standard semi-Lagrangian methods. 

 In essence, the result stated in the next section shows that in this situation, the algorithm \eqref{eq:schema} approximates the exact solution up to an error that is of the order $\mathcal{O}(\delta x + \sqrt{\delta t/N})$. The first term comes from the interpolation, and the second show a variance reduction phenomenon. To obtain this bound, we require an anti-CLF condition as usual for semi-Lagrangian methods. We also assume that $N$ is sufficiently large in some relatively weak sense (see a precise statement below). 

The paper is organized as follows. In Section \ref{SectPres}, we present the method, introduce various notations and state our main result. In Section 3 we give some properties of random matrices arising in a natural way in the definition of the scheme \eqref{eq:schema}, and that are needed in the proof of our main estimate, given in full details Section 4. 
Possible extensions of the method, with for instance Dirichlet boundary conditions, or more complicated PDEs, are evoked in Section \ref{Exten}, together with a few numerical results. In particular, we present simulations for the two-dimensional Burgers equation.

\section{Setting and main result\label{SectPres}}

%For simplicity, we focus on the approximation of a simple equation, which is linear, and associated to a simple stochastic process. Moreover, summability properties will be easily satisfied working in a periodic setting.

We consider the linear heat equation on $(0,1)$, with a smooth initial condition and with periodic boundary conditions; more precisely, we want to approximate the unique periodic solution of the following partial differential equation:
\begin{equation}
\label{eq:hot}
\frac{\partial u}{\partial t}=\frac{1}{2}\frac{\partial^2 u}{\partial x^2},
\end{equation}
such that $u(0,.)=u_0:\R\rightarrow \R$ is a smooth periodic function of period $1$, and such that for any $t\geq 0$ the function $u(t,.)$ is periodic of period $1$. 
For periodic functions, we denote by $L^2$ and $H^1$ the usual function spaces associated with the respective norm and semi-norm
\begin{gather*}
\|f\|_{L^2}^{2}=\int_{0}^{1}|f(x)|^2\dd x,\quad \mbox{and}\quad
|f|_{H^1}^{2}=\int_{0}^{1}|f'(x)|^2\dd x.
\end{gather*}

\subsection{Interpolation operator}

%We need to discretize time and space, and to define an interpolation procedure allowing to recover functions defined on the line from given values at the nodes.

%First, we consider a final time $T>0$, and an integer $M_T$, such that we divide the interval $[0,T]$ into $M_T$ intervals of size $\delta t:=\frac{T}{M_T}$. We are thus interested in approximating the solution $u$ at times $t_n=n\delta t$ for $n\leq M_T$. The constants in the error bounds may depend on the finite time $T$ but not on the discretization time step $\delta t$.

For a given integer $M_S \geq 1$, we discretize the space interval  $(0,1)$ with the introduction of nodes $x_j=j\delta x$ for $j\in S := \left\{0,\ldots,M_S - 1\right\}$, with the condition $x_{M_S}=M_S\delta x=1$. We set $V = \{ (u_j)_{j \in S} \} \subset \R^M$ the set of discrete functions defined on the discrete points of the grid.

%By periodicity, the definition of $x_j$ may be extended for indices in $\mathbb{Z}$. In sums written below, we may forget to mention the set of indices: it must be then considered that we sum over a period: for example we often choose a canonical discrete period $S=\left\{0,\ldots,M_S-1\right\}$, which represents the nodes in $[0,1)$, since by periodicity the information contained in values at index $j=M_S$ already appears for index $j=0$.

We use linear interpolation to reconstruct functions on the whole interval from values at the nodes. We define an appropriate basis made of periodic and piecewise linear functions for $k\in S=\left\{0,\ldots,M_S-1\right\}$. We set for $x \in 
[x_k-1/2,x_k+1/2]$, 
$$
\phi_k(x)=\hat{\phi}(\frac{x-x_k}{\delta x}), \quad \mbox{where}\quad 
\hat{\phi}(x)=\begin{cases}
        0 \text{ if } |x|>1,\\
	1-|x| \text{ if }|x|\leq 1, 
        \end{cases}
$$
and extend the function $\phi_k$ by periodicity on $(0,1)$. 
Hence $\phi_k$ is the unique piecewise linear  periodic function and satisfying
$\phi_k(x_j)= \delta_{kj}$ the Kronecker symbol. Note that we have $\sum_{k\in S}\phi_k(x)=1$ for all $x \in (0,1)$.

%We remark that for any $x\in [0,1]$
%\begin{equation}\label{remphik}
%\phi_k(x)\geq 0 \quad \text{and} \quad \sum_{k\in S}\phi_k(x)=1.
%\end{equation}
%%This only means that the constant function $1$ can be interpolated linearly, and that we recover the initial function.
%
%Moreover we have the following useful property: for any $k\in S$ and any $x\in \R$
%\begin{equation}\label{transphik}
%\phi_{k+1}(x)=\phi_{k}(x-\delta x)
%\end{equation}
%If we denote by $\mathcal{B}_{p,1}(\R)$ the set of all $1$-periodic functions with values in $\R$, and by $\mathcal{C}_{p,1}(\R)$ the set of the continuous $1$-periodic functions, 
We define the following projection and interpolation operators
$$
\mathcal{P}:\begin{cases}
H^1\rightarrow V\\
f\mapsto (f(x_j))_{j \in S}
\end{cases}
\quad \mbox{and}\quad 
\mathcal{I}:\begin{cases}
V \rightarrow H^1\\
u=(u_k)_{k \in S}\mapsto \sum_{k=0}^{M_S-1}u_k\phi_k(x),
\end{cases}
$$
where $H^1$ denote the Sobolev space of periodic function on $(0,1)$. 
Clearly, $\mathcal{P}\circ\mathcal{I}$ is the identity on $V$; nevertheless the distance between the identity and the composition of the operators $\mathcal{I}\circ\mathcal{P}$ depends on the functional spaces and on the norms. Below, we give the estimates that are useful in our setting. We just notice that  $\mathcal{I}\circ\mathcal{P}(\mathds{1})=\mathds{1}$, as $\sum_{k} \phi_k(x) = 1$ for all $x$.

\subsection{Discrete norms}

%@We define appropriate norms on $\R^{M_S}$, which are related to norms for functions defined on continuous intervals. We moreover state some important properties of the operators $\mathcal{P}$ and $\mathcal{I}$ related to the norms, which give estimates of the error made by interpolation.

%Below $J$ denotes any period; for definiteness, we take $J=S:=\left\{0,\ldots,M_S-1\right\}$.
For any $u=(u_j)_{j\in S}$ we define the discrete $L^2$ norm and $H^1$ semi-norm as 
$$
\|u\|_{\ell^2}^{2}=\delta x\sum_{j\in S}u_{j}^{2}
\quad \mbox{and}
\quad 
|u|_{h^1}^{2}=\delta x\sum_{j\in S}\frac{(u_{j+1}-u_{j})^2}{\delta x^2}, 
$$
where we use the extension by periodicity of the sequence $(u_j)$ for the definition of the $h^1$ semi-norm: we thus have $u_{M_S}=u_0$. %As in the continuous case, with $|\hspace{0.1cm}.\hspace{0.1cm}|_{h^1}$ we only define a semi-norm: if $|u|_{h^1}=0$, then there exists $c\in \R$ such that for any $0\leq j\leq M_S-1$ we have $u_j=\frac{1}{M_S}\sum_{k=0}^{M_S-1}u_k$.
We also define a norm with $\|u\|_{h^1}=(\|u\|_{\ell^2}^{2}+|u|_{h^1}^{2})^{1/2}$. 

%Both discrete norm and semi-norm are approximations of the corresponding quantities on the whole interval $[0,1]$:
%\begin{gather*}
%\|f\|_{L^2(0,1)}^{2}=\int_{0}^{1}|f(x)|^2dx, \text{ for }f\in L^{2}(0,1),\\
%|f|_{H^1(0,1)}^{2}=\int_{0}^{1}|f'(x)|^2dx, \text{ for }f\in H^{1}(0,1).
%\end{gather*}
%

With these notations, we have the following approximation results: 
\begin{propo}\label{PropEstimInterp}
There exists a constant $c>0$ such that for any mesh size $\delta x=1/M_S$, and any sequence $u=(u_j) \in V$ we have:
\begin{gather*}
|u|_{h^1}=|\mathcal{I}u|_{H^1}, \quad \mbox{and}\quad
\|u\|_{\ell^2}^{2}=\|\mathcal{I}u\|_{L^2}^{2}+c\delta x^2|u|_{h^1}^{2}.
\end{gather*}
Moreover, for any function $f\in H^1$ we have
$$\|f-(\mathcal{I}\circ\mathcal{P})f\|_{L^2}^{2}\leq c\delta x^2\left(|f|_{H^1}^{2}+|(\mathcal{I}\circ\mathcal{P})f|_{H^1}^{2}\right).$$
\end{propo}

%As a consequence, we also get the following estimate, for any $f\in H^1(0,1)$:
%$$|\|f\|_{L^{2}(0,1)}^{2}-\|\mathcal{P}f\|_{\ell^2}^{2}|\leq c\delta x^2\left(|f|_{H^1(0,1)}^{2}+|(\mathcal{I}\circ\mathcal{P})f|_{H^1(0,1)}^{2}\right).$$

%\begin{rem}
%Thanks to a usual Sobolev embedding, $H^1(0,1)\subset C^0(0,1)$ the set of continuous functions, and therefore $\mathcal{P}f$ is well-defined.% - taking for a function $f$ the unique continuous function $\tilde{f}$ such that $f=\tilde{f}$ almost everywhere.
%\end{rem}

%
\underline{Proof}. 
The first equality follows from a direct computation. 
The second one is proved expanding in the $L^2$ scalar product $\mathcal{I}u=\sum_ku_k\phi_k$, and rewriting the sums in order to make the $h^1$ semi-norm appear: we have
$$
\|\mathcal{I}u\|_{L^{2}}^{2}=\|\sum_{k\in S}u_k\phi_k\|_{L^2}^{2}
=\sum_{k,\ell\in S}u_ku_\ell\left\langle\phi_k,\phi_\ell\right\rangle_{L^2}
=\frac{2\delta x}{3}\sum_{k\in S}u_{k}^{2}+\frac{\delta x}{6}\sum_{k\in S}(u_{k}u_{k+1}+u_{k}u_{k-1}), 
$$
where we define $u_{M_S} = u_0$ and $u_{-1} = u_{M_S - 1}$, that is we extend $u \in V$ by periodicity. We also used the fact that for all $k$, 
$$
\left\langle\phi_k,\phi_\ell\right\rangle_{L^2}=0 \text{ if }\ell\notin\left\{k-1,k,k+1\right\}, \quad 
\left\langle\phi_k,\phi_k\right\rangle_{L^2}=\frac{2\delta x}{3}, \quad \mbox{and} \quad 
\left\langle\phi_k,\phi_{k-1}\right\rangle_{L^2}=\frac{\delta x}{6}.
$$

Now, the equality contains $\|u\|_{h^1}^{2}$ which appears with natural integration by parts - using periodicity:
\begin{align*}
\|u\|_{\ell^2}^{2}-\|\mathcal{I}u\|_{L^2}^{2}&=\frac{\delta x}{6}\sum_{k\in S}\left(u_{k}(u_{k}-u_{k-1})+u_{k}(u_{k}-u_{k+1})\right)\\
&=\frac{\delta x}{6}\sum_{k\in S}\left(u_{k+1}(u_{k+1}-u_{k})+u_{k}(u_{k}-u_{k+1})\right)\\
&=\frac{\delta x}{6}\sum_{k\in S}(u_{k+1}-u_{k})^2=\frac{1}{6}\delta x^2\|u\|_{h^1}^{2}.
\end{align*}

To prove the last estimate, we have for any $j\in S=\left\{0,1,\ldots M_S-1 \right\}$ and for any $x\in[x_{j},x_{j+1}]$, 
\begin{align*}
|f(x)-(\mathcal{I}\circ\mathcal{P}f)(x)|^2&\leq 2(\int_{x_j}^{x}f'(t)\dd t)^2+2(\int_{x_j}^{x}\frac{f(x_{j+1})-f(x_j)}{\delta x}\dd t)^2\\
&\leq 2(x-x_{j})\int_{x_j}^{x_{j+1}}|f'(t)|^2\dd t+2\delta x^2\frac{|f(x_{j+1})-f(x_j)|^2}{\delta x^2},
\end{align*}
using the Cauchy-Schwarz inequality.
Now we integrate over $[x_{j},x_{j+1}]$, and then it remains to take the sum over $j\in S$ of the following quantities:
\begin{equation}\label{IneInterp}
\int_{x_{j}}^{x_{j+1}}|f(x)-(\mathcal{I}\circ\mathcal{P}f)(x)|^2 \dd x\leq \delta x^2\int_{x_j}^{x_{j+1}}|f'(x)|^2\dd x+2\delta x^2\frac{|f(x_{j+1})-f(x_j)|^2}{\delta x^2}\delta x.
\end{equation}
The first term of the right-hand side is controlled with $|f|_{H^1}^{2}$, while the second term involves $|\mathcal{P}f|_{h^1}^{2}=|(\mathcal{I}\circ\mathcal{P})f|_{H^1}^{2}$.
%Now if we take a periodic function $f\in H^{1}(0,1)$, thanks to a Sobolev embedding $f$ is indeed a continuous function, so that $\mathcal{P}(f)$ is well-defined; moreover the previous inequality extends to such functions $f$ by a density argument.
\qed

%We also see that if the function $f$ is of class $\mathcal{C}^1$ on $[0,1]$, we get
%$$|f|_{H^1(0,1)}^{2}+|(\mathcal{I}\circ\mathcal{P})f|_{H^1(0,1)}^{2}\leq \sup_{x\in[0,1]}|f'(t)|^2.$$

%In the sequel, we try to bound the error at time $t_n=n\delta t$ with the $l^2$ norm of discrete quantities, with a bound involving the $h^1$-semi norm of the initial condition.

% discrete norm l^2, discrete semi-norm h^1.
% inequalities

\subsection{Definition of the numerical method}
We consider a final time $T>0$, and an integer $M_T$, such that we divide the interval $[0,T]$ into $M_T$ intervals of size $\delta t:=\frac{T}{M_T}$. We are thus interested in approximating the solution $u(t,x)$ of \eqref{eq:hot} at times $t_n=n\delta t$ and nodes $x_j$ for $n\leq M_T$. %The constants in the error bounds may depend on the finite time $T$ but not on the discretization time step $\delta t$.
In the simple case of Equation \eqref{eq:hot}, for which we have the representation formula $u(t,x)=\E u_0(x+B_t)$, 
the numerical scheme \eqref{eq:schema} is the application $u^n \mapsto u^{n+1}$ from $V$ to itself written 
\begin{equation}
\label{def_unj}
u_{j}^{n+1}=\frac{1}{N}\sum_{m=1}^{N}\left(\sum_{k\in S}u_{k}^{n}\phi_k(x_j+\sqrt{\delta t}\mathcal{N}^{n,m,j})\right),
\end{equation}
where the random variables $\mathcal{N}^{n,m,j}$, indexed by $0\leq n\leq M_T-1$, $1\leq m\leq N$ and $j\in S$ are independent standard normal variables.
More precisely, to avoid an error term due to the approximation of Brownian Motion at discrete times, we require that
\begin{equation}\label{defBMjm}
\sqrt{\delta t}\mathcal{N}^{n,m,j}=B_{(n+1)\delta t}^{(m,j)}-B_{n\delta t}^{(m,j)}
\end{equation}
for some independent Brownian Motions $(B^{(m,j)})$ for $1\leq m\leq N$ and $0\leq j\leq M_S-1$.% This condition does not change the problem.

We start with an initial condition $u^0=(u_{k}^{0}=u_0(x_k))$, which contains the values of the initial condition at the nodes. To obtain simple expressions with products of matrices, we consider that vectors like $u^0$ are column vectors.

We then define the important auxiliary sequence $v^n \in V$ satisfying the following relations:
\begin{equation}\label{def_vnj}
\begin{aligned}
v_{j}^{n+1}&=\frac{1}{N}\sum_{m=1}^{N}\left(\sum_{k\in S}v_{k}^{n}\E[\phi_k(x_j+\sqrt{\delta t}\mathcal{N}^{n,m,j})]\right)\\
&=\sum_{k\in S}v_{k}^{n}\E[\phi_k(x_j+\sqrt{\delta t}\mathcal{N}^{n,m,j})],
\end{aligned}
\end{equation}
with the initial condition $v^0=u^0$.
Indeed, for any $0\leq n\leq M_T$ the vector $v^n$ is the expected value - defined component-wise - of the random vector $u^n$.

%Since we want to work with linear algebra expressions, we consider that vectors like $u^0$ are row vectors.

%Having calculated the approximation at discrete time $n$, we define the values at time $n+1$: for any $j\in S$ the canonical set $S$
%\begin{equation}\label{def_unj}
%u_{j}^{n+1}=\frac{1}{N}\sum_{m=1}^{N}\left(\sum_{k\in S}u_{k}^{n}\phi_k(x_j+\sqrt{\delta t}\mathcal{N}^{n,m,j})\right),
%\end{equation}
%where the random variables $\mathcal{N}^{n,m,j}$, indexed by $0\leq n\leq M_T-1$, $1\leq m\leq N$ and $j\in S$ are independent standard normal variables.
%
%More precisely, to avoid an error term due to the approximation of Brownian Motion at discrete times, we require that
%\begin{equation}\label{defBMjm}
%\sqrt{\delta t}\mathcal{N}^{n,m,j}=B_{(n+1)\delta t}^{(m,j)}-B_{n\delta t}^{(m,j)}
%\end{equation}
%for some independent Brownian Motions $(B^{(m,j)})$ for $1\leq m\leq N$ and $0\leq j\leq M_S-1$.% This condition does not change the problem.

\subsection{Main result} With the previous notations, we have the following error estimate: 

%We also introduce in Definition \ref{DefNormes} a discrete norm $\|.\|_{\ell^2}$ and a discrete semi-norm $|.|_{h^1}$, which are discrete counterparts of the  norm in $L^2(0,1)$ and the semi-norm in $H^1(0,1)$. We consider the full-discretization error, which is decomposed into an error due to randomness - we control its variance in the $l^2$ norm, depending on the $h^1$ semi-norm of the initial condition - and an error only due to deterministic effects - corresponding to the accumulation in time of the interpolation error. The error estimate is given in the following Theorem:

%The Monte-Carlo error is evaluated in the $l^2$ norm, and depends on the $h^1$ semi-norm of the initial condition. 

\begin{theo}\label{TheTheorem}
Assume that the initial condition $u_0$ is of class $\mathcal{C}^2$.
%For any $T>0$, $N_{T}^{0}\in \N$ and $N_{S}^{0}$, there exists a constant $C>0$ such that for any $M_T\geq N_{T}^{0}$ and $N_{S}^{0}$, with $\delta t=\frac{T}{M_T}$ and $\delta x=\frac{1}{M_S}$, and for any $N\geq 1$, the Monte-Carlo error satisfies for any $0\leq n\leq M_T$:
For any $p\in\mathbb{N}$ and any final time $T>0$, there exists a constant $C_p>0$, such that for any $\delta t>0$, $\delta x>0$ and $N\in\mathbb{N}^*$ we have
%\begin{equation}\label{decompeqTh}
%(u(t_n,x_k))_{k}-u^n=(u(t_n,x_k))_{k}-v^n)+(v^n-u^n),
%\end{equation}
%with
\begin{equation}\label{eqTh0}
\sup_{j \in S %j\in \mathbb{N}; 0\leq j\delta x<1
}|u(t_n,x_j)-v_{j}^{n}|\leq C\frac{\delta x^2}{\delta t}\sup_{x\in [0,1]}|u_{0}^{''}(x)|
\end{equation}
and
%$\delta t_0>0$, $\delta x_0$ and $N_0$,
\begin{equation}\label{eqTh}
\E\|u^n-v ^n\|_{\ell^2}^{2}
\leq C_p|u^0|_{h^1}^{2}(1+\frac{\delta x^2}{\delta t})\left(1+\frac{\delta x}{\delta t}+\frac{\delta x^2}{\delta t^2}(1+|\log(\delta t)|)\right)^p\left(\frac{\delta t}{N}+\frac{1}{N^{p+1}} \right).
%\begin{aligned}
%\E\|u^n-v ^n\|_{\ell^2}^{2}&=\delta x\sum_{j; 0\leq j\delta x<1}\E|u_{j}^{n}-v_{j}^{n}|^2\\
%&\leq C_p|u^0|_{h^1}^{2}(1+\frac{\delta x^2}{\delta t})\left(1+\frac{\delta x}{\delta t}+\frac{\delta x^2}{\delta t^2}(1+|\log(\delta t)|)\right)^p\left(\frac{\delta t}{N}+\frac{1}{N^{p+1}} \right).
%\end{aligned}
%\begin{aligned}
%\E\|u^n-v ^n\|_{\ell^2}^{2}&\leq C\|u^0\|_{h^1}^{2}\frac{1+\frac{\delta x^2}{\delta t}}{N}(\delta t+\delta x+\frac{\delta x^2}{\delta t}|\log(\delta t)|)\\
%&+C\|u^0\|_{h^1}^{2}\frac{1+\frac{\delta x^2}{\delta t}}{N^2}(1+\frac{\delta x}{\delta t}+\frac{\delta x^2}{\delta t^2}|\log(\delta t)|).
%\end{aligned}
\end{equation}
\end{theo}

The control of the first part of the error is rather classical, while the estimate on the Monte-Carlo error given by \eqref{eqTh} is more original and requires more attention in its analysis and in its proof.

First, we observe that the estimate is only interesting if a condition of anti-CFL type is satisfied: for some constant $c>0$ we require
$$\frac{\delta x}{\delta t}\text{max}(1,\sqrt{\log(\delta t)})<c.$$
We then identify in \eqref{eqTh} a leading term of size $\frac{\delta t}{N}$, which corresponds to the statistical error in a Monte-Carlo method for random variables of variance $\delta t$, and a remaining term, which goes to $0$ with arbitrary order of convergence with respect to the number of realizations $N$. This second term is obtained via a bootstrap argument. Indeed it is easy to get the classical estimate with $p=0$. The core of the proof is contained in the recursion which allows to increase the order from $p$ to $p+1$; it heavily relies on the spatial structure of the noise and on the choice of the $\ell^2$-norm.

Thanks to \eqref{eqTh} when $p=1$, we see that interpreting Theorem \ref{TheTheorem} as a reduction of variance to a size $\delta t$ is valid: we bound the error with%$$\frac{\delta t}{N}+\frac{1}{N^2}\leq \frac{\delta t^2}{2}+\frac{3}{2N^2}\quad \text{and} \quad
%\text{Var}(\frac{1}{N}\sum_{i=1}^{N}Y_i)=\frac{\text{Var}(Y)}{N}\leq \frac{\text{Var}(Y)^2}{2}+\frac{1}{2N^2}.$$
%We also give a nice interpretation when $p=1$: we have
$$\frac{\delta t}{N}+\frac{1}{N^2}\leq \frac{\delta t^2}{2}+\frac{3}{2N^2},$$
which can be compared with a classical Monte-Carlo bound with the variance $\delta t$: we have for any sample $(Y_1,\ldots,Y_N)$ of a random variable $Y$
$$\text{Var}\left(\frac{1}{N}\sum_{i=1}^{N}Y_i\right)=\frac{\text{Var}(Y)}{N}\leq \frac{\text{Var}(Y)^2}{2}+\frac{1}{2N^2}.$$
%Indeed, the estimate \eqref{eqTh} is better, and shows that the interpretation with a variance $\delta t$ is not exact but quite precise.

The control of the Monte-Carlo error in Theorem \ref{TheTheorem} relies on several arguments. Firstly, the first factor corresponds to the accumulation of the variances appearing at each time step - where two sources of error are identified: the random variables involve a stochastic diffusion process evaluated at time $\delta t$, and an error is introduced by the interpolation procedure. To obtain another factor, we observe that the independence of the random variables appearing for different nodes implies that only diagonal entries of some matrices appear - see \eqref{Lezero}. However, this independence property also complicates the proof: the solutions are badly controlled with respect to the $h^1$ semi-norm. We then propose a decomposition of the error where the number of realizations $N$ appears in the variance with the different orders $1$ and $2$: the first part is controlled by $\delta t$ and $\delta x$, while the second one is only bounded. We finally use recursively this decomposition in order to improve the estimate, with a bootstrap argument.

%More complicated PDEs can be discretized according to the same kind of Monte-Carlo semi-lagrangian method, and we have obtained numerical simulations for various examples: for instance the space variable can be $d$-dimensional, and Dirichlet or Neumann conditions can be imposed on the boundary.

%Extensions of the numerical method with numerical experiments are given in Chapter $5$.

\section{Random matrices}

The definition of the numerical scheme \eqref{def_unj} can be rewritten with matrix notations: for column vectors of size $M_S$ such that $(u^n)_j=u_{j}^{n}$, we see that
\begin{equation}\label{matrixu^n}
u^{n+1}=P^{(n)}u^{n},
\end{equation}
where the entries of square matrix satisfy for any $1\leq j,k\leq M_S$
\begin{equation}\label{defPn}
P^{(n)}_{j,k}=\frac{1}{N}\sum_{m=1}^{N}\phi_k(x_j+\sqrt{\delta t}\mathcal{N}^{n,m,j}).
\end{equation}
Moreover we decompose these matrices into $N$ independent parts: for $1\leq m\leq N$
\begin{equation}\label{defPnm}
P^{(n)}=\frac{1}{N}\sum_{m=1}^{N}P^{(n,m)},
\end{equation}
with the entries $(P^{(n,m)})_{j,k}=\phi_k(x_j+\sqrt{\delta t}\mathcal{N}^{n,m,j})$.

We observe that the matrices $P^{(n,m)}$ are independent; in each one, the rows are independent, however in a row indexed by $j$ two different entries are never independent, since they depend on the same random variable $\mathcal{N}^{n,m,j}$; moreover, the sum of coefficients in a row is $1$.

All matrices $P^{(n,m)}$ have the same law; we define a matrix $Q=\E P^{(n,m)}=\E P^{(n)}$, by taking the expectations of each entry: for any $j,k\in S$
\begin{equation}\label{defQjk}
Q_{j,k}=\E[\phi_k(x_j+\sqrt{\delta t}\mathcal{N}^{n,m,j})].
\end{equation}
The right-hand side above does not depend on $n,m$ since we take expectation. It only depends on $j$ through the position $x_j$, not through the random variable $\mathcal{N}^{n,m,j}$. With these notations, the vectors $v^n$ satisfy the relation $v^{n+1}=  Q v^n$ -  see \eqref{def_vnj} -  and we have for any $n\geq 0$
\begin{equation}\label{schemematrix}
u^n=\prod_{i=0}^{n-1}P^{(i)} u^0=P^{(n-1)}\ldots P^{(0)}u^0, 
\quad \mbox{and} \quad
v^n=Q^nu^0.
\end{equation}

%%%%%%%%%%%%%%%%%%%%%

We only present a few basic properties of the matrices $P^{(n,m)}$, $P^{(n)}$ and $Q$. First, we show that they are stochastic matrices. Second, we control their behavior with respect to the discrete norms and semi-norms. In order to prove the convergence result, we need other more technical properties which are developed during the proof.

%First, all the matrices are stochastic matrices. This is an easy consequence of their definition with the interpolation functions.
\begin{propo}\label{propstosym}
For any $0\leq n\leq M_T-1$ and for any $1\leq m\leq N$, almost surely $P^{(n,m)}$ is a stochastic matrix: for any indices $j,k\in S$ we have $P_{j,k}^{(n,m)}\geq 0$, and for any $j\in S$ we have $\sum_{k\in S}P_{j,k}^{(n,m)}=1$.

For any $0\leq n\leq M_T-1$, $P^{(n)}$ is also a random stochastic matrix.

The matrix $Q$ is stochastic and symmetric - and therefore is bistochastic.
\end{propo}

\underline{Proof}. 
The stochasticity of the random matrices $P^{(n,m)}$ is a simple consequence of their definition \eqref{defPn} and of the relations  $\phi_k(x) \geq 0$ and $\sum_{k \in S} \phi_{kl}(x) = 1$. Since $P^{(n)}$ is a convex sum of the $P^{(n,m)}$, the property for those matrices also holds.

Finally, by taking expectation $Q$ is obviously stochastic; symmetry is a consequence of \eqref{defQjk}, and of the property $\phi_{k+1}(x)=\phi_{k}(x-\delta x)$:
\begin{align*}
Q_{j,k}&=\E[\phi_k(x_j+\sqrt{\delta t}\mathcal{N}^{n,m,j})]
=\E[\phi_0(x_j-x_k+\sqrt{\delta t}\mathcal{N}^{n,m,j})]\\
&=\E[\phi_0(x_k-x_j-\sqrt{\delta t}\mathcal{N}^{n,m,j})]
=\E[\phi_0(x_k-x_j+\sqrt{\delta t}\mathcal{N}^{n,m,k})]=Q_{k,j}, 
%&=\E[\hat{\phi}(\frac{x_j-x_k+\sqrt{\delta t}\mathcal{N}^{n,m,j}}{\delta x})]\\
%&=\E[\hat{\phi}(\frac{x_k-x_j-\sqrt{\delta t}\mathcal{N}^{n,m,j}}{\delta x})]\\
%&=Q_{k,j},
\end{align*}
since $\phi_0$ is an even function, and since the law of $\mathcal{N}^{n,m,j}$ is symmetric and does not depend on $j$. However this symmetry property is not satisfied by the $P$-matrices, because the trajectories of these random variables are different when $j$ changes.
\qed

Thanks to the chain of equalities in the proof above, we see that $Q_{j,k}$ only depends on $k-j$, but we observe that no similar property holds for the matrices $P^{(n,m)}$.

We now focus on the behavior of the matrices with respect to the $\ell^2$-norm. The following proposition is a discrete counterpart of the decreasing of the $L^2$-norm of solutions of the heat equation.
\begin{propo}\label{Pl2}
For any $0\leq n\leq M_T-1$ and for any $1\leq m\leq N$, and for any $u\in V$ we have
$$\E\|P^{(n,m)}u\|_{\ell^2}^{2}\leq \|u\| _{\ell^2}^{2} \quad\text{ and }\quad
%This property is easily extended to the matrices $P^{(n)}$: for any $0\leq n\leq M_T-1$ and for any $u\in \R^{M_S}$ we have
\E\|P^{(n)}u\|_{\ell^2}^{2}\leq \|u\| _{\ell^2}^{2}.$$
\end{propo}

\underline{Proof}. 
According to the definitions above \eqref{matrixu^n} and \eqref{defPnm}, we have for any index $j$
$(P^{(n,m)}u)_j=\sum_{k\in S}P_{j,k}^{(n,m)}u_k$.
Thanks to the previous Proposition \ref{propstosym}, we use the Jensen inequality to get
\begin{align*}
\E\|P^{(n,m)}u\|_{\ell^2}^{2}&=\delta x\sum_{j\in S}\E|(P^{(n,m)}u)_j|^2\\
&\leq \delta x\sum_{j\in S}\sum_{k\in S}\E P_{j,k}^{(n,m)}|u_{k}|^{2}
\leq \delta x\sum_{k\in S}\left(\sum_{j\in S}Q_{j,k}\right)|u_k|^2;
\end{align*}
now we use the properties of the matrix $Q$ - it is a bistochastic matrix according to Proposition \ref{propstosym} - to conclude the proof, since $\sum_{j\in S}Q_{j,k}=1$.
The extension to the matrices $P^{(n)}$ is straightforward.
\qed

The matrix $Q$ satisfies the same decreasing property in the $\ell^2$-norm; moreover we easily obtain a bound relative to the $h^1$-semi norm:
\begin{propo}\label{Ql2}
For any $u\in V$, we have $\|Qu\|_{\ell^2}\leq \|u\| _{\ell^2} \text{  and  } 
|Qu|_{h^1}\leq |u| _{h^1}.$
\end{propo}

\underline{Proof}. 
The proof of the first inequality is similar to the previous situation for the random matrices. To get the second one, it suffices to define a sequence $\tilde{u}$ such that for any $0\leq j\leq M_S-1$ we have $\tilde{u}_j=\frac{u_{j+1}-u_{j}}{\delta x}$ - with the convention $u_{M_S}=u_0$. Then thanks to the properties of $Q$ we have $\widetilde{Qu}=Q\tilde{u}$: for any $j\in S$
\begin{align*}
(\delta x)\widetilde{Qu}_{j}&=(Qu)_{j+1}-(Qu)_{j}=\sum_{k\in S}Q_{j+1,k}u_k-\sum_{k\in S}Q_{j,k}u_k\\
&=\sum_{k\in S}Q_{j,k-1}u_{k}-\sum_{k\in S}Q_{j,k}u_k
=\sum_{k\in S}Q_{j,k}(u_{k+1}-u_{k})=\delta x(Q\tilde{u})_j,
\end{align*}
using a translation of indices with periodic conditions, and the equality $Q_{j+1,k}=Q_{j,k-1}$ as explained above.
As a consequence, we have
$|Qu|_{h^1}=\|\widetilde{Qu}\|_{\ell^2}=\|Q\tilde{u}\|_{\ell^2}\leq \|\tilde{u}\|_{\ell^2}=|u|_{h^1}.$
\qed

It is worth noting that the previous argument can not be used to control $\E|P^{(n,m)}u|_{h^1}$: for a matrix $P=P^{(n,m)}$, the corresponding quantity $\widetilde{Pu}$ can not be easily expressed with $\tilde{u}$. Indeed, given a deterministic $u$, then $(P^{(n,m)}u)_{j}$ and $(P^{(n,m)}u)_{j+1}$ are independent random variables - since they are defined respectively with $\mathcal{N}^{(n,m,j)}$ and $\mathcal{N}^{(n,m,j+1)}$. The only result that can be proved is Proposition \ref{Ph1} below. However, its only role in the sequel is to explain why we can not obtain directly a good error bound; as a consequence, we do not give its proof.
% since it requires Proposition \ref{prop1step}, and since \eqref{one_step_h1} is not used to get Theorem \ref{TheTheorem}, we postpone the proof to Section \ref{sectproofPh1}.
\begin{propo}\label{Ph1}
There exists a constant $C$, such that for any discretization parameters $N\geq 1$, $\delta t=\frac{T}{M_T}$ and $\delta x=\frac{1}{M_S}$, we have for any vector $u\in V$
\begin{equation}\label{one_step_h1}
\E |P^{(0)}u|_{h^1}^{2}\leq (1+C\frac{\delta t+\delta x^2}{N\delta x^2})|u|_{h^1}^{2}.
\end{equation}
\end{propo}

%We can remark that in \eqref{one_step_h1} the control depends on the number of Monte-Carlo simulations: as we guess from the expression of $P^{(0)}u$, the treatment of random variables corresponding to identical or distinct realizations is different. On the one hand, Monte-Carlo independence implies the same behavior as for $Qu$; on the other hand, for a same Monte-Carlo realization independence of positions $j+1$ and $j$ leads to an estimate which does not take into account the smallness of $\delta x$.

Due to independence of matrices involved at different steps of the scheme, the previous inequalities can be used in chain.
% Nevertheless, it is dangerous to use the estimate of Proposition \ref{Ph1}. Indeed, the recursive use of such an inequality may in general lead to a non accurate estimate. For example when $\delta x=\delta t$ we have:
%$$\E |P^{(n-1)}\ldots P^{(0)}u|_{h^1}^{2}\leq  (1+C\frac{\delta t}{N\delta x^2})^n|u|_{h^1}^{2}\leq e^{\frac{CT}{N\delta x^2}}|u|_{h^1}^{2},$$
%for any $n\leq \frac{T}{\Delta t}$. Without any further relation between $N$ and $\delta x$ this estimate is useless.

%Depending on the relations between $\delta x$ and $N$ the right-hand side may explode exponentially fast, and therefore this estimate is not satisfactory. 

We thus observe that the matrices $P^{(k)}$ and $Q$ are quite different, even if $Q=\E P^{(k)}$. On the one hand, the matrix $Q$ is symmetric, and therefore respects the structure of the heat equation - the Laplace operator is also symmetric with respect to the $L^2$-scalar product. On the other hand, the structure of the noise destroys this symmetry for matrices $P^{(k)}$, while it introduces many other properties due to independence - in some sense noise is white in space and implies first that solutions are not regular, but that on the average a better estimate can be obtained.

\section{Proof of Theorem \ref{TheTheorem}}\label{SectProof}

We begin with a detailed proof of \eqref{eqTh}. A proof of the other part of the error \eqref{eqTh0} is given in Section \ref{SectAccu} below. 
%The proof of Theorem \ref{TheTheorem} is based on several arguments. First, we consider the error with respect to the $l^2$-norm at time $n$, so that we have the following expression, for any $0\leq n\leq M_T$:
%The numerical error is controlled with respect to the $l^2$-norm, so that we have the following expression, for any $0\leq n\leq M_T$:
%We recall that the error estimate of \eqref{eqTh} is given in the $l^2$-norm.
Easy computations give the following expression for the part corresponding to the Monte-Carlo error: for any $0\leq n\leq M_T$
\begin{align*}
\delta x\sum_{j=0}^{M_S-1}\text{Var}(u_{j}^{n})&=\delta x\sum_{j=0}^{M_S-1}\E|u_{j}^{n}-v_{j}^{n}|^2
=\E\|u^n-v^n\|_{\ell^2}^{2}
=\delta x\E(u^n-v^n)^{*}(u^n-v^n),
\end{align*}
where the superscript $*$ denotes transposition of matrices.

Since the vectors $u^n$ and $v^n$ satisfy \eqref{schemematrix}, with the same deterministic initial condition $u^0$, we have
\begin{align*}
\E\|u^n-v ^n\|_{\ell^2}^{2}&=\E\|(P^{(n-1)}\ldots P^{(0)}-Q^n)u^0\|_{\ell^2}^{2}\\
&=\delta x (u^0)^{*}\E\left((P^{(n-1)}\ldots P^{(0)}-Q^n)^{*}(P^{(n-1)}\ldots P^{(0)}-Q^n)\right)u^0\\
&=\delta x (u^0)^{*}\E\left((P^{(0)})^{*}\ldots (P^{(n-1)})^{*}P^{(n-1)}\ldots P^{(0)}-(Q^{n})^{*}Q^{n}\right)u^0,
\end{align*}
where the last inequality is a consequence of the relation $\E P^{(k)}=Q$ and of the independence of the matrices $P^{(k)}$.

Therefore we need to study the matrix $S_n=\E\left((P^{(0)})^{*}\ldots (P^{(n-1)})^{*}P^{(n-1)}\ldots P^{(0)}-(Q^{n})^{*}Q^{n}\right)$ given by the expression above, such that
$$\E\|u^n-v ^n\|_{\ell^2}^{2}=\delta x(u^0)^{*}S_nu^0.$$

\subsection{Decompositions of the error}

We propose two decompositions of $S_n$ into sums of $n$ terms, involving products of matrices $P^{(k)}$, of $Q$ and of the difference between two matrices $P^{(k)}$ and $Q$, which corresponds to a one-step error:
\begin{equation}\label{decomp1}
P^{(n-1)}\ldots P^{(0)}-Q^n=\sum_{k=0}^{n-1}P^{(n-1)}\ldots P^{(k+1)}\bigl(P^{(k)}-Q\bigr)Q^k,
\end{equation}
and
\begin{equation}\label{decomp2}
P^{(n-1)}\ldots P^{(0)}-Q^n=\sum_{k=0}^{n-1}Q^{n-1-k}\bigl(P^{(k)}-Q\bigr)P^{(k-1)}\ldots P^{(0)}.
\end{equation}

These decompositions lead to the following expressions for $S_n$ - where we use the independence of the matrices $P^{(k)}$ for different values of $k$:
\begin{align*}
S_n&=\E\sum_{k=0}^{n-1}(Q^{k})^{*}\Bigl(P^{(k)}-Q\Bigr)^{*}(P^{(k+1)})^{*}\ldots (P^{(n-1)})^{*}P^{(n-1)}\ldots P^{(k+1)}\Bigl(P^{(k)}-Q\Bigr)Q^k\\
&=\E\sum_{k=0}^{n-1}(P^{(0)})^{*}\ldots (P^{(k-1)})^{*}\Bigl(P^{(k)}-Q\Bigr)^{*}(Q^{n-1-k})^{*}Q^{n-1-k}\Bigl(P^{(k)}-Q\Bigr)P^{(k-1)}\ldots P^{(0)}
\end{align*}

Therefore we obtain the following expressions for the error:
\begin{equation}\label{decompserr}
\begin{aligned}
\E\|u^n-v ^n\|_{\ell^2}^{2}&=\delta x(u^0)^{*}S_nu^0\\
&=\sum_{k=0}^{n-1}\E\|P^{(n-1)}\ldots P^{(k+1)}\bigl(P^{(k)}-Q\bigr)Q^ku^0\|_{\ell^2}^{2}\\
&=\sum_{k=0}^{n-1}\E\|Q^{n-1-k}\bigl(P^{(k)}-Q\bigr)P^{(k-1)}\ldots P^{(0)}u^0\|_{\ell^2}^{2}.
\end{aligned}
\end{equation}

Before we show how each decomposition is used to obtain a convergence result, we focus on the variance induced by one step of the scheme. In fact, only the second one gives the improved estimate of Theorem \ref{TheTheorem}. Nevertheless, we also get a useful error bound thanks to the first one.

\subsection{One-step variance}

In the previous Section, we have introduced decompositions of the error, and we observed that we need a bound on the error made after each time-step. The following Proposition states that the variance after one step of the scheme is of size $\delta t$ if we consider the $\ell^2$ norm, and that a residual term of size $\delta x^2$ appears due to the interpolation procedure. If we consider $N$ independent realizations, Corollary \ref{cor1step} below states that the variance is divided by $1/N$ if we look at the full matrix of the scheme.
\begin{propo}\label{prop1step}
There exists a constant $C$, such that for any discretization parameters $\delta t=\frac{T}{M_T}$ and $\delta x=\frac{1}{M_S}$, and for any $1\leq m\leq N$ and $0\leq n\leq M_T-1$, we have for any vector $u\in \R^{M_S}$
\begin{equation}\label{1step}
\E\|(P^{(n,m)}-Q)u\|_{\ell^2}^{2}\leq C(\delta t+\delta x^2)|u|_{h^1}^{2}.
\end{equation}
\end{propo}

\begin{cor}\label{cor1step}
For any $0\leq n\leq M_T-1$ and for any vector $u\in \R^{M_S}$, we have
$$\E\|(P^{(n)}-Q)u\|_{\ell^2}^{2}\leq C\frac{(\delta t+\delta x^2)}{N}|u|_{h^1}^{2}.$$
\end{cor}

The proof of the corollary is straightforward, since $P^{(n)}=\frac{1}{N}\sum_{m=1}^{N}P^{(n,m)}$ with independent and identically distributed matrices $P^{(n,m)}$. However, the proof of Proposition \ref{prop1step} is very technical.

%Heuristically, the error is of order $1/2$ with respect to the time step since the left-hand side can be interpreted as the variance of a functional of a diffusion process, evaluated at time $\delta t$. Due to the use of discrete norms and to the interpolation procedure, an error of size $\delta x^2$ also appears.% Finally, the error is controlled with respect to the $h^1$- semi norm of the initial condition $u$, which contains information on the regularity of this vector - and of a function $f$ such that $\mathcal{P}f=u$.

One difficulty of the proof is the dependence of the noise on the position $j$: for different indices $j_1$ and $j_2$, the random variables $(P^{(n,m)}u)_{j_1}$ and $(P^{(n,m)}u)_{j_2}$ are independent. To deal with this problem, for each $j$ we introduce an appropriate auxiliary function and we analyze the error on each interval $[x_{j},x_{j+1}]$ separately. We also need to take care of some regularity properties of the functions - they are $H^1$ functions, piecewise linear, but they are not in general of class $\mathcal{C}^1$ - in order to obtain bounds involving the $h^1$ and $H^1$ semi-norms.

\vspace{0.5cm}

\underline{Proof of Proposition \ref{prop1step}}. 
To simplify the notations, we assume that $n=0$ and that $m=1$%; using \eqref{defBMjm}, we could write explicitly similar expressions with the true indices, but the final result remains the same. Therefore
so that we only work with one matrix $P$ with entries
$$P_{j,k}=\phi_k(x_j+B_{\delta t}^{j}),$$
where the $B^j$ are independent Brownian Motions.

We define the following auxiliary periodic functions: for any $x\in \R$
\begin{equation}\label{defVx}
V(x)=\E \II u(x+B_{\delta t}^{j}),
\end{equation}
and for any index $0\leq j\leq M_S-1$
\begin{equation}\label{defUjx}
U^{(j)}(x)=\II u(x+B_{\delta t}^{j}).
\end{equation}

We observe that since we take expectation in \eqref{defVx} the index $j$ plays no role there. Moreover we have the following relations for any $j\in S$:
\begin{gather*}
V(x_j)=(Qu)_j\quad \mbox{and}\quad 
U^{(j)}(x_j)=(Pu)_j,\quad \mbox{but}\quad 
U^{(j)}(x_{j+1})\neq (Pu)_{j+1}.
\end{gather*}
The last relation is the reason why we need to introduce different auxiliary functions $U^{(j)}$ for each index $j$.

We finally introduce the following function depending on two variables: for any $0\leq t\leq \delta t$ and $x\in \R$,
\begin{equation}\label{defVtx}
\mathcal{V}(t,x)=\E \II u(x+B_{t}),
\end{equation}
for some standard Brownian Motion $B$. This function is solution of the backward Kolmogorov equation associated with the Brownian Motion, with the initial condition $\mathcal{V}(0,.)=\II u$, and for $t>0$
$$\partial_t \mathcal{V}=\frac{1}{2}\partial_{xx}^{2}\mathcal{V}.$$
Moreover we have $\mathcal{V}(\delta t,.)=V$.% Indeed, $\mathcal{V}$ is solution of the heat equation in a periodic setting, which we are trying to discretize with a semi-lagrangian method, but with a different initial condition obtained by interpolation.

We have the following expression for the mean-square error, integrated over an interval $[x_{j},x_{j+1}]$: for any index $j\in S$
\begin{equation}\label{ErrorH1}
\int_{x_j}^{x_{j+1}}\E|U^{(j)}(x)-V(x)|^2\dd x=\int_{0}^{\delta t}\int_{x_j}^{x_{j+1}}\E|\partial_x\mathcal{V}(\delta t-s,x+B_{s}^{j})|^2\dd x\dd s.
\end{equation}
The proof of this identity is as follows. First, thanks to smoothing properties of the heat semi-group, for any $t>0$ the function $\mathcal{V}(t,.)$ is smooth. Using It\^o formula, with the Brownian Motion $B^{(j)}$ corresponding to the function $U^{(j)}$,
$$\dd \mathcal{V}(\delta t-s,x+B_{s}^{j})=\partial_x \mathcal{V}(\delta t-s,x+B_{s}^{j})\dd B_{s}^{j},$$
for $0\leq s\leq \delta t-\epsilon$ and for any $\epsilon\in(0,\delta t)$, and the isometry property implies
$$\E|\mathcal{V}(\delta t,x)-\mathcal{V}(\epsilon ,x+B_{\delta t-\epsilon}^{j})|^2=\int_{0}^{\delta t-\epsilon}|\partial_x \mathcal{V}(\delta t-s,x+B_{s}^{j})|^2\dd s.$$
We integrate over $x\in[x_{j},x_{j+1}]$, and we then pass to the limit $\epsilon\rightarrow 0$, since $\mathcal{V}(0,.)=\II u$ is a piecewise linear function. Moreover, we use the identity $\mathcal{V}(\delta t,.)=V$. We observe that in the right-hand side of the last equality we take expectation, so that we replace $B^{j}$ with the Brownian Motion $B$, which does not depend on $j$.

%More precisely, assuming that $\psi$ is a smooth function, if we define $\Psi(t,x)=\E \psi(x+B_t)$, then $\Psi$ is solution of the Kolmogorov equation $\partial_t\Psi=\frac{1}{2}\partial_{xx}^{2}\Psi$, with the initial condition $\Psi(0,.)=\psi$. Thanks to It\^o formula, we obtain the identity $d\Psi(\delta t-s,x+B_s)=\partial_x\Psi(\delta t-s,x+B_s)dB_s$ for $0\leq s\leq \delta t$, and by the isometry property $\E|\psi(x+B_{\delta t})-\Psi(\delta t,x)^2=\int_{0}^{\delta t}|\partial_x\Psi(\delta t-s,x+B_s)|^2ds$. To conclude, we integrate over $[x_{j},x_{j+1}]$, and observe that the identity can be extended to the case $\phi=\II u\in H^1$. We also recall that due to smoothing properties of the heat semi-group, for any $t>0$ the function $\mathcal{V}(t,.)$ is smooth.

%Due to smoothing properties of the heat semi-group, for any $t>0$ the function $\mathcal{V}(t,.)$ is smooth. Moreover, as a consequence of It\^o formula, we can obtain the following expression for the mean-square error: for any index $j\in S$
%$$\int_{x_j}^{x_{j+1}}\E|U^{(j)}(x)-V(x)|^2dx=\int_{0}^{\delta t}\int_{x_j}^{x_{j+1}}\E|\partial_x\mathcal{V}(\delta t-s,x+B_{s}^{j})|^2dsdx.$$
%We have to integrate over a space interval $[x_{j},x_{j+1}]$ to have a well-defined right-hand side - since the initial condition $V(0,.)=\II u$ is not smooth but only $H^1$. The estimate can be obtained first for smooth initial conditions - and without integration - and then for more general ones thanks to a regularization and passing to the limit.

%We also notice that the right-hand side of \eqref{ErrorH1} does not depend on $j$.
Summing over indices $j\in S$, we then get, thanks to an affine change of variables $y=x+B_s$
\begin{align*}
\sum_{j\in S}\int_{x_j}^{x_{j+1}}\E|U^{(j)}(x)-V(x)|^2\dd x&=\int_{0}^{\delta t}\int_{0}^{1}\E|\partial_x\mathcal{V}(\delta t-s,x+B_{s})|^2 \dd x\dd s\\
&=\int_{0}^{\delta t}\int_{0}^{1}|\partial_x\mathcal{V}(\delta t-s,x)|^2\dd x\dd s=\int_{0}^{\delta t}|\mathcal{V}(\delta t-s,.)|_{H^1}^{2}\dd s\\
&\leq\int_{0}^{\delta t}|\mathcal{V}(0,.)|_{H^1}^{2}\dd s=\delta t|\II u|_{H^1}^{2}=\delta t|u|_{h^1}^{2}.
\end{align*}
%We used periodicity of the functions - with an affine change of variables $y=x+B_s$ in the space integral - and the decreasing of the $H^1$ semi-norm for solutions of the heat equation with periodic boundary conditions on $(0,1)$.

%Now we claim that for some constant $C$ - which does not depend on the parameters or on $u$ - we have
%\begin{equation}\label{claim0}
%\Big|\E\|Pu-Qu\|_{\ell^2}^{2}-\sum_{j\in S}\int_{x_j}^{x_{j+1}}\E|U^{(j)}(x)-V(x)|^2dx\Big|\leq C\delta x^2|u|_{h^1}^{2}.
%\end{equation}
%The proof of \eqref{claim0} is done in two steps: first we show that 
The inequality \eqref{1step} is then a consequence of the two following estimates: first,
\begin{equation}\label{pro1}
\sum_{j\in S}\int_{x_j}^{x_{j+1}}\E\Big|U^{(j)}(x)-V(x)-\II\circ\PPP(U^{(j)}-V)(x)\Big|^2\dd x\leq C\delta x^2|u|_{h^1}^{2},
\end{equation}
and second we show that
\begin{equation}\label{pro2}
\Big|\E\|Pu-Qu\|_{\ell^2}^{2}-\sum_{j\in S}\int_{x_j}^{x_{j+1}}\E|\II\circ\PPP(U^{(j)}-V)(x)|^2\dd x\Big|\leq C\delta x^2|u|_{h^1}^{2}.
\end{equation}

%Each function  $u^{(j)}$ is almost surely a $H^1$ function.
%For each realization of the random variable $B_{\delta t}^{j}$, the function $U^{(j)}$ belongs to $H^1(0,1)$. We use the inequality 
%\eqref{IneInterp} from the proof of Proposition \ref{PropEstimInterp} on each interval $[x_j,x_{j+1}]$ and we have
To get \eqref{pro1}, we use the inequality \eqref{IneInterp} on each interval $[x_j,x_{j+1}]$, for a fixed realization of $B_{\delta t}^{j}$:
\begin{align*}
\int_{x_j}^{x_{j+1}}|U^{(j)}(x)-V(x)&-\II\circ\PPP(U^{(j)}-V)(x)|^2\dd x\leq C\delta x^2\int_{x_{j}}^{x_{j+1}}|\partial_x(U^{(j)}-V)(x)|^2\dd x\\
&+C\delta x\delta x^2\frac{|[U^{(j)}(x_{j+1})-V(x_{j+1})]-[U^{(j)}(x_{j})-V(x_{j})]|^2}{\delta x^2}.
\end{align*}

%We treat separately each term in the right-hand side above.
Taking the sum over indices $j\in S$ and expectation, we see that
%\begin{multline*}
\begin{align*}
\sum_{j\in S}\int_{x_{j}}^{x_{j+1}}\E|\partial_x(U^{(j)}-V)(x)|^2\dd x&\leq
2\sum_{j\in S}\int_{x_{j}}^{x_{j+1}}\E|\partial_x(\II u)(x+B_{\delta t}^{j})|^2\dd x
+\sum_{j\in S}\int_{x_j}^{x_{j+1}}|\partial_xV(x)|^2\dd x\\
&\leq 2(|\II u|_{H^1}^{2}+|\mathcal{V}(\delta t,.)|_{H^1}^{2})
\leq 4|\II u|_{H^1}^{2}=4|u|_{h^1}^{2},
\end{align*}
%\end{multline*}
since $V=\mathcal{V}(\delta t,.)$. Indeed, taking expectation allows to consider a single Brownian Motion $B$, without $j$-dependence. % and using the fact that taking expectation implies that we can take a Brownian motion $B$ which does not depend on $j$, and then by an obvious change of variables we obtain the result.

%We treat the other part in the same way:
We now decompose the remaining term as follows:
\begin{align*}
\frac{|[U^{(j)}(x_{j+1})-V(x_{j+1})]-[U^{(j)}(x_{j})-V(x_{j})]|^2}{\delta x^2}&\leq 2\frac{|U^{(j)}(x_{j+1})-U^{(j)}(x_j)|^2}{\delta x^2}\\
&+2\frac{|V(x_{j+1})-V(x_{j})|^2}{\delta x^2}.
\end{align*}
With the second part, using Proposition \ref{Ql2} we see that
%\begin{align*}
$$\delta x\sum_{j\in S}\frac{|V(x_{j+1})-V(x_{j})|^2}{\delta x^2}=\delta x\sum_{j\in S}\frac{|(Qu)_{j+1}-(Qu)_{j}|^2}{\delta x^2}
=|Qu|_{h^1}^{2}
\leq |u|_{h^1}^{2}.$$
%\end{align*}
To treat the first part, we make the fundamental observation that for a fixed $j\in S$, the same noise process $B^{j}$ is used to compute all values $U^{(j)}(x)$ when $x$ varies. As a consequence, we can use a pathwise, almost sure version of the argument leading to the proof of Proposition \ref{Ql2} which concerns the behavior of $Q$ with respect to the $h^1$ semi norm.% We recall that the obstruction for a similar result on matrices $P$ is spatial roughness of the noise process, due to independence of the Gaussian random variables - see Proposition \ref{Ph1}.
%The first part can be treated in the same way: for a fixed $j\in S$, the same noise process $B^{j}$ is used to compute all values $U^{(j)}(x)$ when $x$ varies. Therefore the argument leading to the estimate of the behavior of $Q$ with respect to the $h^1$ semi norm can be adapted, whereas it does not lead to an estimate for $P$. More precisely, we have
\begin{align*}
U^{(j)}(x_{j+1})-U^{(j)}(x_j)&=\sum_{k\in S}u_k[\phi_k(x_{j+1}+B_{\delta t}^{j})-\phi_k(x_{j}+B_{\delta t}^{j})]\\
&=\sum_{k\in S}u_k[\phi_{k-1}(x_{j}+B_{\delta t}^{j})-\phi_k(x_{j}+B_{\delta t}^{j})]\\
&=\sum_{k\in S}[u_{k+1}-u_{k}]\phi_{k}(x_j+B_{\delta t}^{j}),
\end{align*}
using the relation $\phi_{k+1}(x)=\phi_{k}(x-\delta x)$ and an integration by parts.

Now summing over indices $j\in S$ and using the Jensen inequality - thanks to Proposition \ref{propstosym} - we obtain
\begin{align*}
\delta x\sum_{j\in S}\E\frac{|U^{(j)}(x_{j+1})-U^{(j)}(x_j)|^2}{\delta x^2}&\leq \delta x\sum_{k\in S,j\in S}\E\phi_{k}(x_j+B_{\delta t}^{j})\frac{|u_{k+1}-u_{k}|^2}{\delta x^2}\\
&\leq \delta x\sum_{k\in S}\frac{|u_{k+1}-u_{k}|^2}{\delta x^2}=|u|_{h^1}^{2}.
\end{align*}

Having proved \eqref{pro1}, we now focus on \eqref{pro2}. We have, since $\PPP(U^{(j)}-V)_j=[(P-Q)u]_j$
\begin{align*}
|\sum_{j\in S}\int_{x_{j}}^{x_{j+1}}|\II\circ\PPP (U^{(j)}-V)(x)|^2dx&-\delta x\sum_{j}|[(P-Q)u]_j|^2|\\
&\leq C\delta x^2\delta x\sum_{j\in S}\frac{|(U^{(j)}-V)(x_{j+1})-(U^{(j)}-V)(x_{j})|^2}{\delta x^2}.
\end{align*}
It remains to take expectation and to conclude like for \eqref{pro1}.

\qed

\subsection{Proof of Theorem \ref{TheTheorem}}

As we have explained in the introduction, we consider that $\delta x$ is controlled by $\delta t$ thanks to a anti-CFL condition.  Roughly, from Proposition \ref{prop1step} we thus see that the variance obtained after one step of the scheme is of size $\delta t$, and that the error depends on the solution through the $h^1$ semi-norm. Moreover, from Propositions \ref{Ql2} and \ref{Ph1} we remark that the behaviors of the matrices $Q$ and $P^{(n)}$ with respect to this semi-norm are quite different.

Using the first decomposition of the error in \eqref{decompserr}, we use in chain the bounds given above in Propositions \ref{Pl2}, \ref{Ql2} and \ref{prop1step} and Corollary \ref{cor1step}: 
\begin{align*}
\E\|u^n-v ^n\|_{\ell^2}^{2}&=\sum_{k=0}^{n-1}\E\|P^{(n-1)}\ldots P^{(k+1)}(P^{(k)}-Q)Q^ku^0\|_{\ell^2}^{2}\leq \sum_{k=0}^{n-1}\E\|(P^{(k)}-Q)Q^ku^0\|_{\ell^2}^{2}\\
&\leq \sum_{k=0}^{n-1}C\frac{(\delta t+\delta x^2)}{N}|Q^k u^0|_{h^1}^{2}\leq \sum_{k=0}^{n-1}C\frac{(\delta t+\delta x^2)}{N}|u^0|_{h^1}^{2}\\
&\leq C\frac{1+\delta x^2/\delta t}{N}|u^0|_{h^1}^{2}.
\end{align*}

If the continuous problem is initialized with the function $u_0$, which is periodic and of class $\mathcal{C}^1$, then $u^0=\mathcal{P}u_0$ satisfies $|u^0|_{h^1}\leq \sup_{x\in[0,1]}|u_{0}'(x)|$. Moreover we assume that an anti-CFL condition is satisfied, so that the term $\delta x^2/\delta t$ is bounded. As a consequence, we find a classical Monte-Carlo estimate, where the error does not decrease when $\delta t$ goes to $0$ and is only controlled with the number of realizations:
\begin{equation}\label{resultDecomp1}
\E\|u^n-v ^n\|_{\ell^2}^{2}\leq C\frac{1+\delta x^2/\delta t}{N}|u^0|_{h^1}^{2}.
\end{equation}

In fact, \eqref{resultDecomp1} shows that the variances obtained at each time step can be summed to obtain some control of the variance at the final time. To get an improved bound, we thus need other arguments.

The main observation is that using independence of rows in the $P$-matrices, we only need to focus on diagonal terms
$$\sup_{j\in S}\left((Q^\ell)^{*}Q^\ell\right)_{jj}=\sup_{j\in S}\left(Q^{2\ell}\right)_{jj},$$
for indices $0\leq \ell=n-k-1\leq n-1$. We recall that indeed $Q$ is a symmetric matrix, so that $(Q^\ell)^{*}Q^\ell=Q^{2\ell}$.

More precisely, the error can be written
\begin{align*}
\E\|u^n-v ^n\|_{\ell^2}^{2}&=\delta x (u^0)^{*}S_nu^0=\delta x\sum_{i,j\in S}u_{i}^{0}(S_n)_{i,j}u_{j}^{0}\\
&=\delta x\sum_{k=0}^{n-1}\sum_{i,j\in S}u_{i}^{0}\E\bigl((A_k)^{*}(P^{(k)}-Q)Q^{2(n-1-k)}(P^{(k)}-Q)A_k\bigr)_{i,j}u_{j}^{0}, 
\end{align*}
where for simplicity we use the notation $A_k:=P^{(k-1)}\ldots P^{(0)}$.
We compute for any $i,j\in S$, using the independence properties at different steps
\begin{align*}
\E((A_k)^{*}&(P^{(k)}-Q)Q^{2(n-1-k)}(P^{(k)}-Q)A_k)_{i,j}\\
&=\sum_{k_1,k_2,k_3,k_4\in S}\E[(A_k)_{k_1,i}(P^{(k)}-Q)_{k_2,k_1}(Q^{2(n-1-k)})_{k_2,k_3}(P^{(k)}-Q)_{k_3,k_4}(A_k)_{k_4,j}]\\
&=\sum_{k_1,k_2,k_3,k_4\in S}\E[(A_k)_{k_1,i}(A_k)_{k_4,j}]\E[(P^{(k)}-Q)_{k_2,k_1}(P^{(k)}-Q)_{k_3,k_4}](Q^{2(n-1-k)})_{k_2,k_3}.
\end{align*}
The observation is now that if $k_2\neq k_3$, then the independence of the random variables for different nodes implies that
\begin{equation}\label{Lezero}
\E[(P^{(k)}-Q)_{k_2,k_1}(P^{(k)}-Q)_{k_3,k_4}]=0,
\end{equation}
since it is the covariance of two independent random variables - see \eqref{defPn}. Moreover, when $k_2=k_3$ we see that $((Q^{(n-1-k)})^{*}Q^{(n-1-k)})_{k_2,k_3}$ only depends on $n-k-1$, due to invariance properties of the equation. Therefore we rewrite the former expansion in the following way:
\begin{equation}\label{EstimMagic}
\begin{aligned}
\E\|&u^n-v ^n\|_{\ell^2}^{2}=\delta x\sum_{k=0}^{n-1}\sum_{i,j\in S}u_{i}^{0}\E\bigl((A_k)^{*}(P^{(k)}-Q)Q^{2(n-1-k)}(P^{(k)}-Q)A_k\bigr)_{i,j}u_{j}^{0}\\
&=\delta x\sum_{k=0}^{n-1}\sum_{i,j\in S}u_{i}^{0}u_{j}^{0}(Q^{2(n-1-k)})_{1,1}\sum_{k_2\in S}\E\left[\left((P^{(k)}-Q)A_k\right)_{k_2,i}\left((P^{(k)}-Q)A_k\right)_{k_2,j}\right]\\
&=\sum_{k=0}^{n-1}\left(Q^{2(n-1-k)}\right)_{1,1}\E\|(P^{(k)}-Q)P^{(k-1)}\ldots P^{(0)}u^0\|_{\ell^2}^{2}.
\end{aligned}
\end{equation}

%$$\E\|u^n-v ^n\|_{\ell^2}^{2}\leq \sum_{k=0}^{n-1}\sup_{j\in S}\left((Q^{n-1-k})^{*}Q^{n-1-k}\right)_{jj}\E\|(P^{(k)}-Q)P^{(k-1)}\ldots P^{(0)}u^0\|_{\ell^2}^{2}.$$

We thus have to control $(Q^{2\ell})_{1,1}=(Q^{2\ell})_{j,j}$ for any $j\in S$. The following Lemma \ref{TheLemma} gives a control of this expression. The first estimate means that the coefficients $Q_{j_1,j_2}^{2\ell}$ are approximations of the solution of the PDE at time $2\ell\delta t$, at position $j_2$, starting from the initial condition $\phi_{j_1}$, with an error due to interpolation. The second estimate is fundamental in the proof of the Theorem, since it allows to introduce an additional factor $\delta x$; however, we need to treat carefully the denominator.

% Since the matrices $P^{(n,m)}$ are not symmetric, such an analysis is only possible for the second decomposition, and gives an essential estimate.
%We can think of $(Q^{2l})_{jj}$ as an approximation by the semi-lagrangian method of the value of the solution of the partial differential equation, at time $2l\delta t\leq 2T$, at position $j$, and starting from the initial condition $\phi_j$. More precisely, we can prove the following result.

\begin{lemme}\label{TheLemma}
There exists a constant $C$ such that for any discretization parameters $\delta t=\frac{T}{M_T}$  and $\delta x=\frac{1}{M_S}$, we have for any $1\leq \ell\leq M_T-1$ and for any $0\leq j_1,j_2\leq M_S-1$
\begin{equation}\label{eqLemma1}
|Q^{2\ell}_{j_1,j_2}-\E\phi_{j_1}(x_{j_2}+B_{2\ell\delta t})|\leq C\frac{\delta x^2}{\delta t}(1+|\log(\delta t)|).
\end{equation}
Moreover, for any $j\in S$, we have for any $1\leq \ell\leq M_T$
\begin{equation}\label{eqLemma2}
\E\phi_{j}(x_{j}+B_{2\ell\delta t})\leq C\frac{\delta x}{\sqrt{2\ell\delta t}}.
\end{equation}
\end{lemme}

\begin{rem}%The estimate \eqref{eqLemma2} only holds for $\ell\geq 1$, since for $\ell=0$ we have $\phi_j(x_j)=1$. This case is exceptional and we can consider that the common behavior is given by the Lemma.

%In the first estimate \eqref{eqLemma1}, the singularity when $\delta t\rightarrow 0$ with a fixed $\delta x$ is linked to the necessity of proving a uniform bound to control the interpolation error, and to the use of regularization properties of the heat equation: the test functions used here as initial conditions $\phi_j$ are bounded but non smooth, while the solution for positive times becomes smooth.

The singularities when $\delta t\rightarrow 0$ with a fixed $\delta x$ come from the use of regularization of the heat semi-group - when we consider the $\phi_j$'s as initial conditions.

For the second estimate \eqref{eqLemma2}, we make two important remarks. First, the constant $C$ depends on the final time $T$, and we cannot directly let $\ell$ tend to $+\infty$: we have
$$\lim_{\ell\rightarrow +\infty}\E\phi_{j}(x_{j}+B_{2\ell\delta t})=\int_{x_j-1/2}^{x_j+1/2}\phi_j(x)\dd x=\delta x\neq 0.$$
Second, from \eqref{eqLemma2} we get for any $\ell>0$ and for any fixed $\delta t$
$$\lim_{\delta x\rightarrow 0}\E\phi_{j}(x_{j}+B_{2\ell\delta t})=0,$$
while we know that for a fixed $\delta x>0$ and a fixed $\ell$, we have
$$\lim_{\delta t\rightarrow 0}\E\phi_{j}(x_{j}+B_{2\ell\delta t})=\phi_j(x_j)=1.$$
These two behaviors are different and from \eqref{eqLemma2} we see the kind of relations that the parameters $\delta x$ and $\delta t$ must satisfy for obtaining one convergence or the other.
\qed
\end{rem}
%Moreover, from the theory of weak approximations, as explained in the introduction, we could choose other random variables than the Gaussian ones $\mathcal{N}^{n,m,j}$: for example uniform random variables $\mathcal{U}^{n,m,j}$ on the interval $[-1,1]$ or on the two-point set $\left\{-1,1\right\}$. The expression $\E\phi_{j}(x_{j}+B_{2l\delta t})$ should be replaced by $\E\phi_{j}(x_{j}+\sqrt{\delta t}S_{2l})$, for some random walk $(S_l)$, with independent and uniformly distributed increments. Then the above limit would be $0$ if and only if $\PP(S_{2l}=0)=0$. Therefore the choice of random variable with a density, like the Gaussian ones, is essential to get a bound like \eqref{eqLemma2}.

\underline{Proof of Lemma \ref{TheLemma}}. 
For any $0\leq \ell\leq 2M_T$, we define
$$M_\ell=\sup_{i,j\in S}|(Q^{\ell})_{i,j}-\E\phi_j(x_i+B_{\ell\delta t})|,$$
where $(B_t)_{t\geq 0}$ is a standard Brownian Motion.

We have $M_0=0$, and by definition of $Q$ we also have $M_1=0$.

We define some auxiliary functions $W_j$, for any index $j\in S$: for any $x\in\R$ and any $t\geq 0$
$$W_j(t,x)=\E\phi_j(x+B_t).$$
$W_j$ is solution of the heat equation, with periodic boundary conditions and initial condition $\phi_j$. For any $t>0$, $W_j(t,.)$ is therefore a smooth function - thanks to regularization properties of the heat semi-group - and since $\phi_j$ is bounded by $1$ we easily see that we have the following estimates, for some constant $C$:
\begin{equation}\label{regWj}
\begin{gathered}
\|\partial_xW_j(t,.)\|_{\infty}\leq \frac{C}{\sqrt{t}}\quad \mbox{and}\quad 
\|\partial_{xx}^{2}W_j(t,.)\|_{\infty}\leq \frac{C}{t}.
\end{gathered}
\end{equation}

We now prove the following estimate on the sequence $(M_\ell)$: for any $1\leq \ell\leq M_T-1$
\begin{equation}\label{recMk}
M_{\ell+1}\leq M_{\ell}+C\frac{\delta x^2}{\ell\delta t}.
\end{equation}
The error comes from the interpolation procedure which is made at each time step.

For any $i,j\in S$, Markov property implies that
\begin{align*}
(Q^{\ell+1})_{i,j}-\E\phi_j(x_i+B_{(\ell+1)\delta t})&=\sum_{k\in S}Q_{i,k}(Q^{\ell})_{k,j}-\E W_j(\ell\delta t,x_i+B_{\delta t})\\
&=\sum_{k\in S}Q_{i,k}(Q^{\ell})_{k,j}-\E\II\circ\PPP(W_j(\ell\delta t,.))(x_i+B_{\delta t})\\
&+\E[\II\circ\PPP(W_j(\ell\delta t,.))-W_j(\ell\delta t,.)](x_i+B_{\delta t}).
\end{align*}
For the first term, we remark that it is bounded by $M_\ell$; indeed we see that
\begin{align*}
\E\II\circ\PPP(W_j(\ell\delta t,.))(x_i+B_{\delta t})&=\sum_{k\in S}W_j(\ell\delta t,x_k)\E\phi_k(x_i+B_{\delta t})\\
&=\sum_{k\in S}Q_{i,k}\E\phi_j(x_k+B_{\ell\delta t}).
\end{align*}
To conclude, it remains to use the stochasticity of the matrix $Q$: entries are positive, and their sum over each line is equal to $1$.

The second term is bounded using the following argument:
$$\|\II\circ\PPP(W_j(\ell\delta t,.))-W_j(\ell\delta t,.)]\|_{\infty}\leq C\delta x^2\|\partial_{xx}^{2}W_j(\ell\delta t,.)\|_{\infty}\leq C\frac{\delta x^2}{\ell\delta t},$$
according to well-known interpolation estimates and to \eqref{regWj}.

From \eqref{recMk}, using $M_1=0$ we obtain for any $1\leq \ell\leq M_T$
\begin{align*}
M_\ell&\leq C\frac{\delta x^2}{\delta t}\sum_{k=1}^{M_T-1}\frac{1}{k}\leq C\frac{\delta x^2}{\delta t}(|\log(T)|+|\log(\delta t)|).
\end{align*}
which gives the result, with a constant depending on $T$.

Now we prove the second estimate of the Lemma. Thanks to the relation $\phi_{k+1}(x)=\phi_{k}(x-\delta x)$  we see that the left-hand side does not depend on $j\in S$; moreover we expand the calculation of the expectation using the periodicity of the function $\phi_j$ and relation definition of $\hat\phi$, the description of its support as $x_j+\bigcup_{k\in \mathbb{Z}}[k-\delta x,k+\delta x]$: we get for $1\leq \ell\leq M_T$
\begin{align*}
\E\phi_{j}(x_{j}+B_{2\ell\delta t})&=\sum_{k\in \mathbb{Z}}\frac{1}{\sqrt{2\pi \ell\delta t}}\int_{k-\delta x}^{k+\delta x}\hat{\phi}(\frac{z-k}{\delta x})e^{-|z|^2/(2\ell\delta t)}\dd z\\
&\leq \frac{1}{\sqrt{2\pi \ell\delta t}}\sum_{k\in \mathbb{Z}}\int_{k-\delta x}^{k+\delta x}e^{-|z|^2/(2\ell\delta t)}\dd z\\
&\leq \frac{1}{\sqrt{2\pi \ell\delta t}}\sum_{k\in \mathbb{Z}}\int_{k-\delta x}^{k+\delta x}e^{-|z|^2/(2T)}\dd z\\
&\leq \frac{1}{\sqrt{2\pi \ell\delta t}}\sum_{k\in \mathbb{Z}}C\delta xe^{-k^2/(2T)}\\
&\leq C\frac{\delta x}{\sqrt{2\ell\delta t}}.
\end{align*}
\qed

The estimate of Lemma \ref{TheLemma} is now used in \eqref{EstimMagic}, and we obtain:
\begin{equation}\label{defEronde}
\begin{aligned}
\sum_{k=0}^{n-1}\left(Q^{2(n-1-k)}\right)_{1,1}&\leq C\sum_{k=0}^{n-2}\left(\frac{\delta x^2}{\delta t}(1+|\log(\delta t)|)+\frac{\delta x}{\sqrt{(n-1-k)\delta t}}\right)+1\\
&\leq C\left(\frac{\delta x^2}{\delta t^2}(1+|\log(\delta t)|)+\frac{\delta x}{\delta t}+1\right)=:C\mathcal{A}.
\end{aligned}
\end{equation}

To conclude one more argument is necessary: we need to apply Proposition \ref{prop1step} in order to sum the variances. However this involves the quantity $\E|P^{(k-1)}\ldots P^{(0)}u^0|_{h^1}^{2}$, which is badly controlled according to Proposition \ref{Ph1}: for example, when $\delta t=\delta x$ the accumulation only implies that
\begin{align*}
\E |P^{(k-1)}\ldots P^{(0)}u|_{h^1}^{2}&\leq  (1+C\frac{\delta t}{N\delta x^2})^{k}|u|_{h^1}^{2}\leq e^{\frac{CT}{N\delta x^2}}|u|_{h^1}^{2},
\end{align*}
for any $k\leq \frac{T}{\delta t}$. We recall that this bad behavior of the matrices $P^{(n)}$ with respect to the $h^1$-semi norm is a consequence of the independence of the random variables for different nodes, whereas this independence property is essential to get the improved estimate, since it allows to use the second estimate of Lemma \ref{TheLemma}.

\begin{rem}\label{remMatK}%\red{Laisser ou enlever?}\black
Instead of considering Gaussian random variables which are independent with respect to the spatial index $j$, we could more generally introduce - like in \cite{JLL} - a correlation matrix $K$, and try to minimize the variance with respect to the choice of $K$. Here we have chosen $K$ as the identity matrix, so that the noise is white in space; the error bound \eqref{eqTh} we obtain is a nontrivial consequence of an averaging effect due to this choice - see \eqref{Lezero}. A natural question - which is not answered here - would be to analyze the situation for general $K$: do we still improve the variance, and can we get more regular solutions?
\end{rem}

The solution we propose relies on the following idea: if above we could replace $P^{(k-1)}\ldots P^{(0)}$ with $Q^k$, we could easily conclude. Another error term appears, which is controlled by $1/N$ instead of $1/\sqrt{N}$. More precisely, independence properties yield for $k\geq 1$
\begin{equation}\label{decomp_ess}
\begin{aligned}
\E\|(P^{(k)}-Q)P^{(k-1)}\ldots P^{(0)}u^0\|_{\ell^2}^{2}&=\E\|(P^{(k)}-Q)Q^ku^0\|_{\ell^2}^{2}\\
&+\E\|(P^{(k)}-Q)\left(P^{(k-1)}\ldots P^{(0)}-Q^k\right)u^0\|_{\ell^2}^{2}.
\end{aligned}
\end{equation}
The roles of the different terms are as follows. On the one hand, the first term gives the part of size $\frac{\delta t}{N}$, thanks to Lemma \ref{TheLemma}: according to Corollary \ref{cor1step} and to Proposition \ref{Ql2}, we have for any $k\geq 1$ with $k\delta t\leq T$
\begin{equation}
\begin{aligned}
\E\|(P^{(k)}-Q)Q^ku^0\|_{\ell^2}^{2}&\leq C\frac{\delta t+\delta x^2}{N}|Q^ku^0|_{h^1}^{2}\leq C\frac{\delta t+\delta x^2}{N}|u^0|_{h^1}^{2}.
\end{aligned}
\end{equation}

On the other hand, the second term is now used to improve recursively the error estimate, since we have
\begin{equation}
\E\|(P^{(k)}-Q)\left(P^{(k-1)}\ldots P^{(0)}-Q^k\right)u^0\|_{\ell^2}^{2}\leq \frac{C}{N}\E\|\left(P^{(k-1)}\ldots P^{(0)}-Q^k\right)u^0\|_{\ell^2}^{2}.
\end{equation}
The independence of realizations at step $k$ gives the factor $\frac{1}{N}$; we remark that we cannot use the estimation of the one-step variance given by Corollary \ref{cor1step}: otherwise we would need to control $\E\|\left(P^{(k-1)}\ldots P^{(0)}-Q^k\right)u^0\|_{h^1}^{2}$.

Using also \eqref{defEronde} and \eqref{decomp_ess} into \eqref{EstimMagic}, we see that
\begin{equation}\label{inequ_bootstrap}
\begin{aligned}
\sup_{n\in\mathbb{N},n\delta t\leq T}\E\|u^n-v^n\|_{\ell^2}^{2}&\leq C\delta t\mathcal{A}\frac{1+\frac{\delta x^2}{\delta t}}{N}|u^0|_{h^1}^{2}
+C\frac{\mathcal{A}}{N}\sup_{n\in\mathbb{N},n\delta t\leq T}\E\|u^n-v^n\|_{\ell^2}^{2}.
\end{aligned}
\end{equation}
The proof of the Theorem now reduces to the study of the following recursive inequalities, for $p\geq 0$
$$E^{(p+1)}\leq C\delta t\mathcal{A}\frac{1+\frac{\delta x^2}{\delta t}}{N}|u^0|_{h^1}^{2}+\frac{C\mathcal{A}}{N}E^{(p)},$$
with an initialization $E^{(0)}=C\frac{\mathcal{B}}{N}$, according to \eqref{resultDecomp1}, with the notation $\mathcal{B}:=(1+\frac{\delta x^2}{\delta t})|u^0|_{h^1}^{2}$. We remark that the control of the matrices $P^{(k)}$ and $Q$ with respect to the $l^2$-norm leads to another possibility for the initialization: $E^{(0)}=2\|u^0\|_{\ell^2}^{2}$; we observe that the recursion then yields the same kind of estimate.

We finally easily prove that for any $p\geq 0$ there exists a constant $C_p\geq 1$ such that
\begin{equation}\label{finalEp}
\sup_{n\in\mathbb{N},n\delta t\leq T}\E\|u^n-v^n\|_{\ell^2}^{2}\leq C_p\left(\frac{\mathcal{A}^{p}\mathcal{B}}{N^{p+1}}+\mathcal{A}\mathcal{B}\frac{\delta t}{N}\right),
\end{equation}
and the proof of Theorem \ref{TheTheorem} is finished.

%We now show how we can control each term. In order to use the result of Lemma \ref{TheLemma}, we need to treat the cases $k=n-1$ and $k<n-1$ in the sum apart.

\begin{rem}
If we consider the equation
$\frac{\partial u}{\partial t}=\frac{\nu}{2}\frac{\partial^2 u}{\partial x^2}$
with a viscosity parameter $\nu>0$, the quantities $\mathcal{A}$ and $\mathcal{B}$ appearing in the proof are transformed into
\begin{gather*}
\mathcal{A}_{\nu}=(1+\frac{\delta x}{\sqrt{\nu}\delta t}+\frac{\delta x^2}{\nu\delta t^2}(1+|\log(\delta t)|))\quad \mbox{and}\quad 
\mathcal{B}_{\nu}=(\nu+\frac{\delta x^2}{\delta t})|u^0|_{h^1}^{2}.
\end{gather*}

% the Monte-Carlo estimates \eqref{eqTh} and \eqref{resultDecomp1} are transformed into
%\begin{align*}
%\E\|u^n-v ^n\|_{\ell^2}^{2}&=\delta x\sum_{j; 0\leq j\delta x<1}\E|u_{j}^{n}-v_{j}^{n}|^2\\
%&\leq C|u^0|_{h^1}^{2}(\nu+\frac{\delta x^2}{\delta t})\left(1+\frac{\delta x}{\sqrt{\nu}\delta t}+\frac{\delta x^2}{\nu\delta t^2}(1+|\log(\delta t)|)\right)\left(\frac{\delta t}{N}+\frac{1}{N^2} \right).
%\end{align*}
%and
%$$\E\|u^n-v ^n\|_{\ell^2}^{2}\leq C\frac{\nu+\delta x^2/\delta t}{N}|u^0|_{h^1}^{2},$$
where the constant $C$ does not depend on $\nu$.

The first change in the proof concerns the analysis of the one-step variance: in \eqref{1step}, the right-hand side is replaced by $C(\nu\delta t+\delta x^2)$. We observe that the error due to interpolation remains the same.

The second change concerns Lemma \ref{TheLemma}, where we use some regularization properties thanks to gaussian noise: when $\nu$ goes to $0$  the estimates degenerates.

As a consequence, we may observe that the estimate \eqref{eqTh} gives a bound valid for a fixed value of $\nu$, while \eqref{resultDecomp1} becomes more interesting  when $\nu$ is small compared with the discretization parameters.
\end{rem}

\subsection{Accumulation of the interpolation error}\label{SectAccu}

To obtain Theorem \ref{TheTheorem}, it remains to control the deterministic part of the error of the scheme, without the discretization of the expectation with the Monte-Carlo method. We thus need to prove \eqref{eqTh0}:

%The conclusion of Theorem \ref{TheTheorem} concerns the variance in the numerical scheme: under appropriate conditions on the parameters $\delta t$ and $\delta x$, the Monte-Carlo error decreases when $\delta t$ and $\delta x$ tend to $0$. However, the comparison of the numerical and the exact solutions involves other sources of error, as explained in the Introduction: the interpolation and the SDE discretization induce error terms which need to be controlled. For instance, we propose a uniform estimate:

for any $n\in \mathbb{N}$ such that $n\delta t\leq T$, and for any $j\in \mathbb{N}$ with $0\leq x_j=j\delta x<1$, we have
\begin{equation}\label{error_sansMC}
|u(n\delta t,x_j)-v_{j}^{n}|\leq C\frac{\delta x^2}{\delta t}\sup_{x\in [0,1]}|u_{0}^{''}(x)|,
\end{equation}
where $u$ is the exact solution and where $v^n$ is defined by \eqref{def_vnj}.% In other words, this estimate corresponds to the error in the theoretical situation when we assume that expectations are computed exactly.

%In more general situations, another term would appear, linked to the weak error in the discretization of the SDE: we would then bound the error with $\delta t^\beta+\frac{\delta x^\gamma}{\delta t}$, where $\beta$ is the weak order of the numerical scheme, and $\gamma$ is the order of approximation by interpolation in the corresponding norm.

Since $\|u(n\delta t,x_{.})-v^n\|_{\ell^2}\leq \sup_{j}|u(n\delta t,x_{j})-v_{j}^{n}|$, we easily obtain an estimate in the $l^2$-norm. Therefore, the conditions imposed on $\delta x$ and $\delta t$ by \eqref{eqTh} are not restrictive, and can be seen as consequences of the semi-lagrangian framework.

The proof of \eqref{error_sansMC} in our context is as follows: using the exact representation formula and its discrete counterpart \eqref{def_vnj}, we have
\begin{align*}
u((n+1)\delta t,x_j)-v_{j}^{n+1}&=\E u(n\delta t,x_j+B_{\delta t})-\E \sum_{k\in S}v_{k}^{n}\phi_k(x_j+B_{\delta t})\\
&=\sum_{k\in S}(u(n\delta t,x_k)-v_{k}^{n})\E\phi_k(x_j+B_{\delta t})\\
&+\E\left(u(n\delta t,x_j+B_{\delta t})-\sum_{k\in S}u(n\delta t,x_k)\phi_k(x_j+B_{\delta t})\right),
\end{align*}
where $B_{\delta t}$ is a Brownian Motion at time $\delta t$.  

It is easy to see that
$$|\sum_{k\in S}(u(n\delta t,x_k)-v_{k}^{n})\E\phi_k(x_j+B_{\delta t})|\leq \sup_{k\in S}|u(n\delta t,x_k)-v_{k}^{n}|,$$
and we see that the other term depends on the interpolation error:
\begin{align*}
|\E[u(n\delta t,x_j+B_{\delta t})-\sum_{k\in S}u(n\delta t,x_k)&\phi_k(x_j+B_{\delta t})]|\leq \sup_{x\in[0,1]}|u(n\delta t,x)-\II\circ \mathcal{P} u(n\delta t,.)(x)|\\
&\leq C\delta x^2\sup_{x\in[0,1]}|\frac{\partial^2 u}{\partial x^2}(n\delta t,x)|\leq C\delta x^2\sup_{x\in[0,1]}|\frac{\partial^2 u}{\partial x^2}(0,x)|.
\end{align*}

To conclude, we remark that for $n=0$ we have $u(0,x_j)=v_{j}^{0}$.

%\begin{rem}\label{remEuler}
%In more general situations, for instance when the function $c$ is non zero, another error term should be added to \eqref{eqTh0}. It corresponds to the use of a time discretization method on the SDE- for instance of Euler type.
%\end{rem}

\section{Numerical results and extensions}\label{Exten}

\subsection{Illustration of Theorem \ref{TheTheorem}}
\begin{figure}[!h]
\begin{center}
\rotatebox{0}{\resizebox{!}{0.3\linewidth}{
   %P1D0
   %\includegraphics{figures_finales/periodic/P1D0/P1D0.pdf}}}
   %P1D0bis
   %\includegraphics{figures_finales/periodic/P1D0/P1D0bis.pdf}}}
   %P1D0ter
   \includegraphics{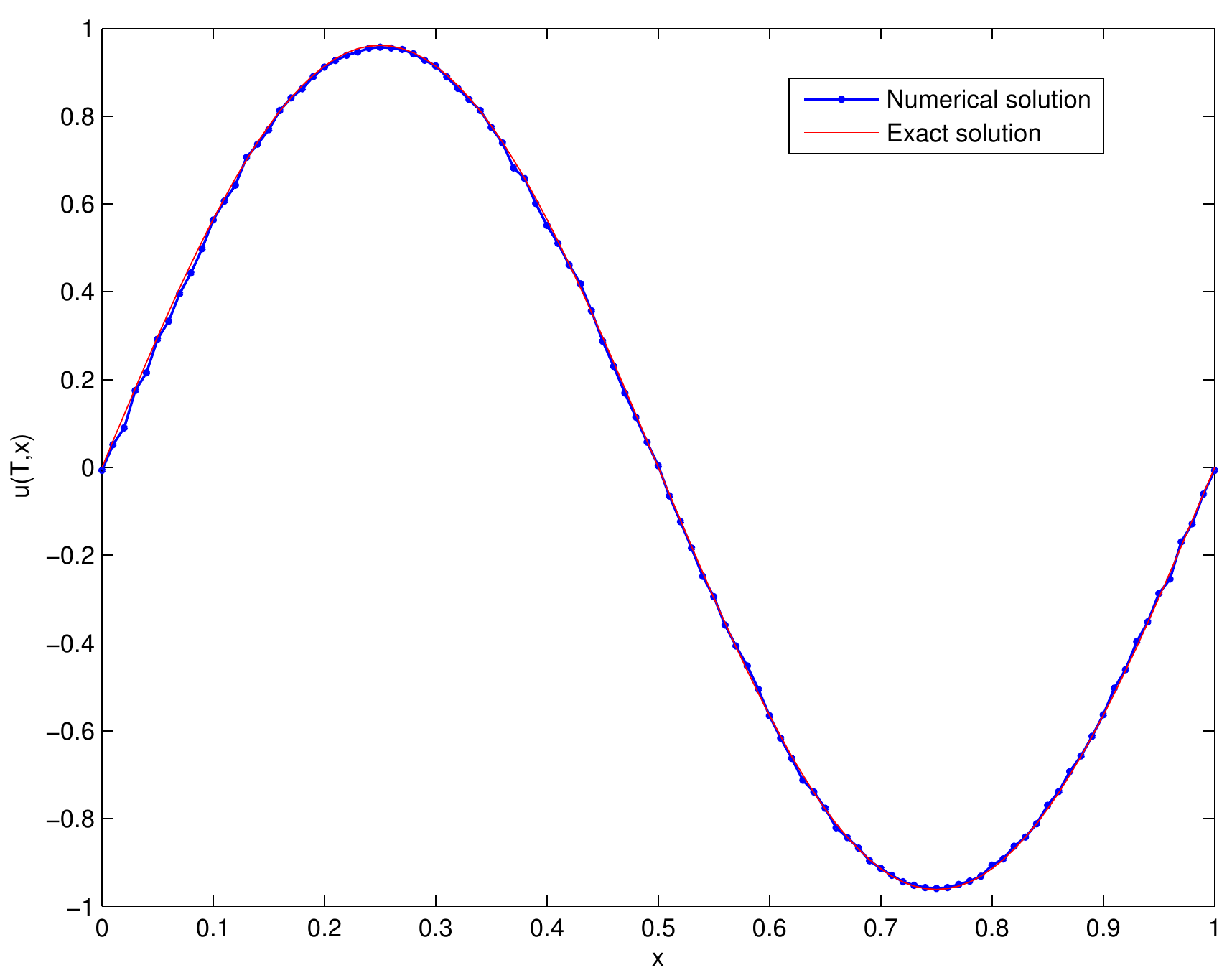}}}
\end{center}
\caption{Solution at time $T=0.1$ with $\delta t=\delta x=1/N=0.01$}
\setlength\belowcaptionskip{0.01ex}
\label{P1D0}
\end{figure}

The first numerical example we consider is a simulation of the solution of the heat equation in the spatial domain $(0,1)$ in periodic setting. We introduce the viscosity parameter $\nu$ so that the problem is
\begin{equation}
\label{eq:num}
\frac{\partial u(t,x)}{\partial t}={\nu}\frac{\partial^2 u(,x)}{\partial x^2}, \text{ for } t>0,x\in (0,1),\quad 
u(0,x)=u_0(x) \text{ for }x\in (0,1), 
\end{equation}
with the boundary condition 
$u(t,1)=u(t,0)$ for $t\geq 0$. 
For the numerical simulation of Figure \ref{P1D0}, we choose $\nu=0.01$, and $u_0(x)=\sin(2\pi x)$. The exact solution satisfies $u(t,x)=\exp(-4\pi^2\nu t)\sin(2\pi x)$. The discretization parameters are
$\delta t=\delta x=0.01$ and 
$N=100$. 

\begin{figure}[!h]
\begin{center}
\rotatebox{0}{\resizebox{!}{0.4\linewidth}{
   \includegraphics{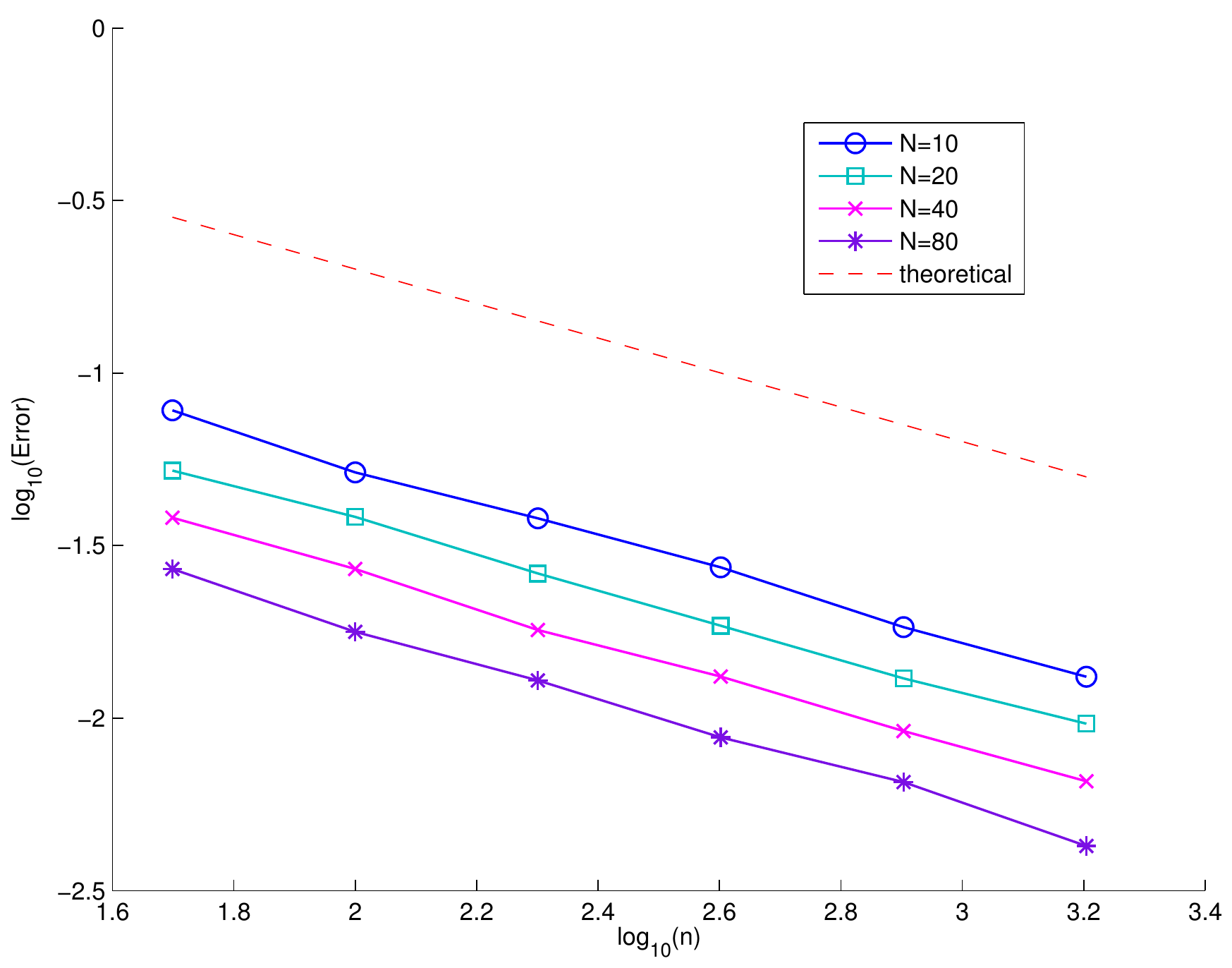}}}
\end{center}
\caption{Error for periodic boundary conditions when $\delta t=\delta x=1/n$, in logarithmic scales.}
\setlength\belowcaptionskip{0.01ex}
\label{errperio}
\end{figure}

The bound of Theorem \ref{TheTheorem} is illustrated with Figure \ref{errperio}, where we represent the error in logarithmic scales for different values of the parameters.

We study the convergence of the scheme, with a numerical simulation which confirms the order of convergence with respect to the parameters $\delta t=\delta x$ of the Monte-Carlo error. The final time is $T=0.1$, the viscosity is $\nu=0.1$ and the initial condition is $u_0(x)=\cos(2\pi x)$. We compare the numerical solution $u^n$ with the exact solution; we only observe the Monte-Carlo error, which is dominant with respect to the deterministic part of the error according to Theorem \ref{TheTheorem}. The mean-square error in the $\ell^2$ norm is estimated with a sample of size $20$.

The error in Figure \ref{errperio} is represented in logarithmic scales. The parameters $\delta t$ and $\delta x$ are equal and satisfy $\delta t=\delta x=\frac{1}{n}$ for the following values $n=50,100,200,400,800,1600,3200$. Each line is obtained when we draw the logarithm of the Error as a function of $\log_{10}(n)$, for a fixed value of $N\in\left\{10,20,40,80\right\}$. The dot-line represents a straight-line with slope $-1/2$.

%In this experiment, the initial condition is $u_0(x)=\cos(2\pi x)$, and the viscosity is $\nu=1$. The error is observed at time $T=0.1$. The expectation $\E\|u^n-v^n\|_{\ell^2}$ is evaluated by a Monte-Carlo method with $400$ realizations.

%We consider two choices for the size of the Monte-Carlo discretization: $N=10$ and $N=20$. The time and space discretization parameters satisfy $\delta t=\delta x=\frac{1}{n}$, with the values $n=200,400,800,1600$.

%We observe an error of size $\sqrt{\delta t}{N}$ in this situation.

This experiment confirms that the Monte-Carlo error is of order $1/2$ with respect to the parameters when $\delta t=\delta x$, as \eqref{eqTh} claims. Indeed, the shift between the lines when $N$ varies also corresponds to the size $1/\sqrt{N}$ of the Monte-Carlo error.

%
%
%In the Introduction, we have described the principle of the hybrid semi-lagrangian method for the equation \eqref{defEDP}, on the domain $\R^d$. We have then studied the case where the coefficient $c$ is null; moreover we have focused on the case of periodic boundary conditions. To go further, we would like to present how the method can be modified in order to consider two examples of new situations: first, when we replace periodic with homogeneous Dirichlet boundary conditions; second, when we approximate some non-linear equations, for instance of Burgers type.
%
%For these two situations, we do not provide a theoretical analysis of the method; we only give some numerical simulations in the next Section \ref{SectNum}. The derivation of error bounds and of possibly more competitive methods will be the content of future works.
%
%

\subsection{The method for Dirichlet boundary conditions}

%The value of the solution on the boundary $\partial D$ - assumed to be of class $\mathcal{C}^2$ - of the bounded domain can also be imposed, given the so-called Dirichlet boundary conditions. Thus the equation is
We would like now show how it is possible to adapt our method in the case of Dirichlet boundary conditions. Let us consider the equation \eqref{eq:num}, but with boundary conditions $u(t,x)=0$ for $t> 0$ and $x\in\partial D= \{ 0,1\}$. 
%We now consider the PDE in a domain $D$ with Dirichlet boundary conditions:
%\begin{gather*}
%\partial_t u=\mathcal{L}u\\
%u(t,.)=u_0\\
%u(t,x)=0, \text{ for } t> 0 \text{ and } x\in\partial D.
%\end{gather*}
The representation formula then involves the family of the first-exit times of the process $X_t^x = x + \sqrt{\nu} B_t$ starting from the different points of the domain: If we define
$\tau^{x}=\inf\left\{t>0;X_t^x\in D^c\right\}$, then the solution satisfies
%$$u(t,x)=\E u_0(X_{t\wedge \tau^x}^{x})=\E \left[u_0(X_{t}^{x})1_{t\leq \tau^{x}} \right];$$
\begin{equation}
\label{eq:stopping}
u(t,x)=\E \left[u_0(X_{t}^{x})\mathds{1}_{t\leq \tau^{x}} \right];
\end{equation}
the stochastic process is killed when it reaches the boundary. Note that this formula extends to more general PDE of the form \eqref{defPDE0} with the associated process \eqref{defEDS0}. 

%If we replace the previous equation by $\partial_t u=\mathcal{L}u+f$ for some function $f$ depending only in the space variable, we still get a representation formula:
%$$u(t,x)=\E [u_0(X_{t}^{x})\mathds{1}_{t< \tau^{x}}+\int_{0}^{t\wedge\tau^x}f(X_{s}^{x})ds].$$

The numerical approximation becomes more complicated, since we also need an accurate approximation of the stopping times. This problem is well-known, and solutions have been proposed in \cite{gobet} and \cite{mann} for the computation of \eqref{eq:stopping} at a given point $x$ using time discretizations of the stochastic process $X_t^x$. 
%In the following paragraphs, we describe elements for the construction of an algorithm.
%The main point of the method focuses on the following situation: it is possible that at the times $n\delta t$ and $(n+1)\delta t$ the simulated values $X_n$ and $X_{n+1}$ both belong to the domain $D$, however an exit of the domain has been possible in continuous time between $n\delta t$ and $(n+1)\delta t$. An additional test is introduced in the scheme to take into account this situation, so that convergence is improved.% This test is described more precisely below for the one-dimensional case; it is based on the knowledge on the law of exit of the diffusion process between $n\delta t$ and $(n+1)\delta t$, conditionally to its values at these times - see formula\eqref{TBB}.

%Nevertheless, in practice this technique requires that the time step $\delta t$ is small; such a restriction is not compatible with the anti-CFL conditions which are naturally associated with semi-lagrangian methods. We propose to refine the time intervals $[n\delta t,(n+1)\delta t]$ in order to improve the approximation of the processes together with the associated stopping times. The interpolation procedure is not concerned with this refinement, and only appears at times $n\delta t$.

In our case, we take advantage of the semi-lagrangian context to do a refinement near the boundary: for a discretization between the times $t_n$ and $t_{n+1}$, we introduce a decomposition of the domain into an "interior" zone and a "boundary" zone, with different treatments. In the boundary zone, we refine in time and use a subdivision of $[n\delta t,(n+1)\delta t]$ of mesh size $\tau\leq \delta t$ and we use a possibly different value $N_b$ for the number of Monte-Carlo realizations. Moreover, following \cite{gobet} and \cite{mann}, we introduce an exit test in the boundary zone, based on the knowledge of the law of exit of the diffusion process. 

In the interior part, less care is necessary and we can take $\tau=\delta t$ and $N_i<N_b$ for the size of the sample. 
%The number of boundary points, the refined time step $\tau$, and the sample sizes $N_b$ and $N_i$ are specified for each simulation.

\begin{figure}[!h]
\begin{center}
\rotatebox{0}{\resizebox{!}{0.4\linewidth}{
   \includegraphics{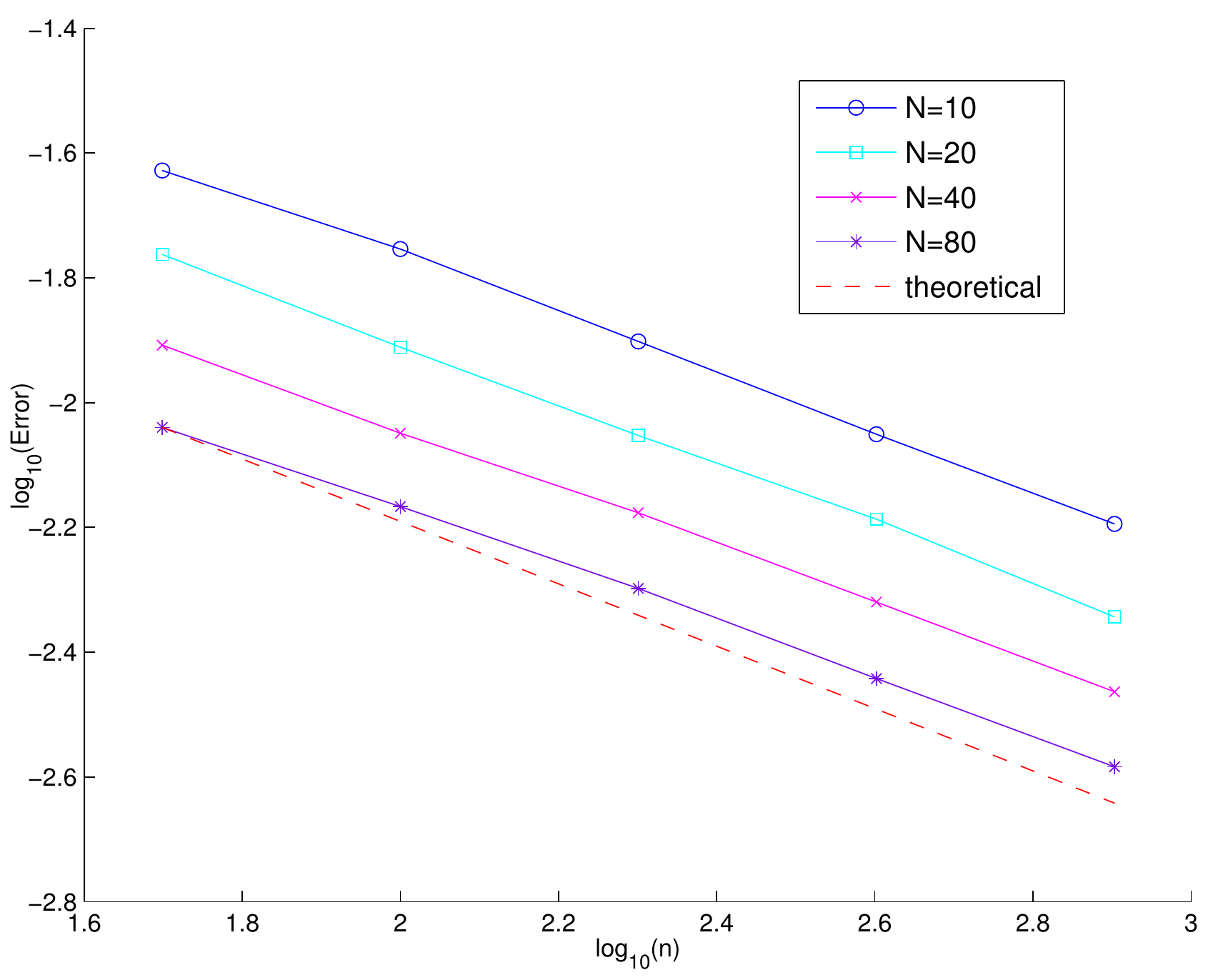}}}
\end{center}
\caption{Error for Dirichlet boundary conditions when $\delta t=\delta x=1/n$, in logarithmic scales.}
\setlength\belowcaptionskip{0.01ex}
\label{err_dir}
\end{figure}

We give in Figure \ref{err_dir} the result of investigations on the convergence of the method when Dirichlet boundary conditions are applied. We draw in logarithmic scales the error in terms of $n=1/\delta t=1/\delta x$, with $n=50,100,200,400,800$, with different values of the Monte-Carlo parameter $N_i=10,20,40,80$. We have chosen on the interval $(-1,1)$ the initial function $u_0(x)=\sin(\pi\frac{x+1}{2})$, with the viscosity $\nu=0.1$. The boundary zone is made of the intervals $(-1,-0.9)$ and $(0.9,1)$, where we take $\tau = \delta t/10$ and $N_b = 10 N_i$. The solutions are computed until time $T=0.1$. Like in the case of periodic boundary conditions, the statistical error is dominant with respect to the other error terms; we compare with the exact solution, and to estimate the variance we use a sample of size $100$. 
%It is important to notice that we have considered $SUB=10$ and that $N_b=10N_i=10N$.

The observation of Figure \ref{err_dir} shows that the Monte-Carlo error depends on the parameter $\delta t=\delta x$; the comparison with the "theoretical" line with slope $-1/2$ indicates a conjecture that the error is also of order $1/2$, like for the periodic case. The shift between the curves for different values of $N$ corresponds in the error to a factor $1/\sqrt{N}$.

\subsection{The method for some non-linear PDEs}

We present a simple method to obtain approximations of the solution of the viscous Burgers equation in dimension $d=2$
$$\frac{\partial u}{\partial t}+(u.\nabla) u=\nu\Delta u+f.$$
It is defined on the domain $(-1,1)^2$, with homogeneous Dirichlet boundary conditions - periodic ones would also have been possible. Compared with the situations described so far, we add a forcing term $f$, which may depend on time $t$, position $x$ and the solution $u$. 

%We must face a new difficulty due to the nonlinear term, which does not allow to construct directly a diffusion process. The strategy is to consider that during a time step the transport occurs with a constant velocity computed with the previous step, so that we obtain a linear parabolic equation on each subinterval, which admits a probabilistic representation formula.

As explained in the Introduction, we construct approximations $u^n$ of the solution at discrete times $n\delta t$, introducing functions $v^n$ such that for any $n\geq 0$ with the following semi-implicit scheme:
\begin{equation}\label{discBurgers}
\frac{\partial v^{n+1}}{\partial t}+(u^n.\nabla) v^{n+1}=\nu\Delta v^{n+1}+f^n,
\end{equation}
for any time $n\delta t\leq t\leq (n+1)\delta t$ and $x\in D$. The initial condition is $v^{n+1}(n\delta t,.)=u^{n}=v^{n}(n\delta t,.)$.
The discrete-time approximation then satisfies $u^0=u_0$ and $u^{n}=v^{n}(n\delta t,.)$. The forcing term here satisfies $f^n(t,x)=f(n\delta t,x,u^n(x))$.

On each subinterval $[n\delta t,(n+1)\delta t]$, we have
$$v^{n+1}(t,x)=\E[v^{n+1}(n\delta t,X_{t}^{x})\mathds{1}_{t<\tau^x}+\int_{n\delta t}^{t\wedge\tau^x}f^n(X_{s}^{x})ds],$$
where the diffusion process $X$ satisfies
\begin{gather*}
dX_{t}^{x}=-u^n(X_{t}^{x})dt+\sqrt{2\nu}dB_t,
X_{n\delta t}^{x}=x.
\end{gather*}
The stopping times $\tau^{x}$ represent the first exit time of the process in the time interval $[n\delta t,(n+1)\delta t]$. Since $v^{n+1}(n\delta t,.)=u^n$, the scheme only requires the knowledge of the approximations $u^n$.

For the numerical simulations, we take the initial condition to be $0$, and the forcing is $f(t,x)=(-\sin(\pi t)\sin(\pi x)\sin(\pi y)^2,-\sin(\pi t)\sin(\pi x)^2\sin(\pi y))$. The viscosity parameter is $\nu=0.001$. The time step satisfies $\delta t=0.02$, and the spatial mesh size is $\delta x=0.04$. The "interior" zone is $(-0.8,+0.8)^2$, where $N_i=10$; on the "boundary" zone, we have $N_b=100=10N_i$, and $\tau=0.002=\delta t/10$.

\begin{figure}[h!]
\centering
\subfigure[$u_1,t=0.5$]{%
     \includegraphics[width=0.35\linewidth]{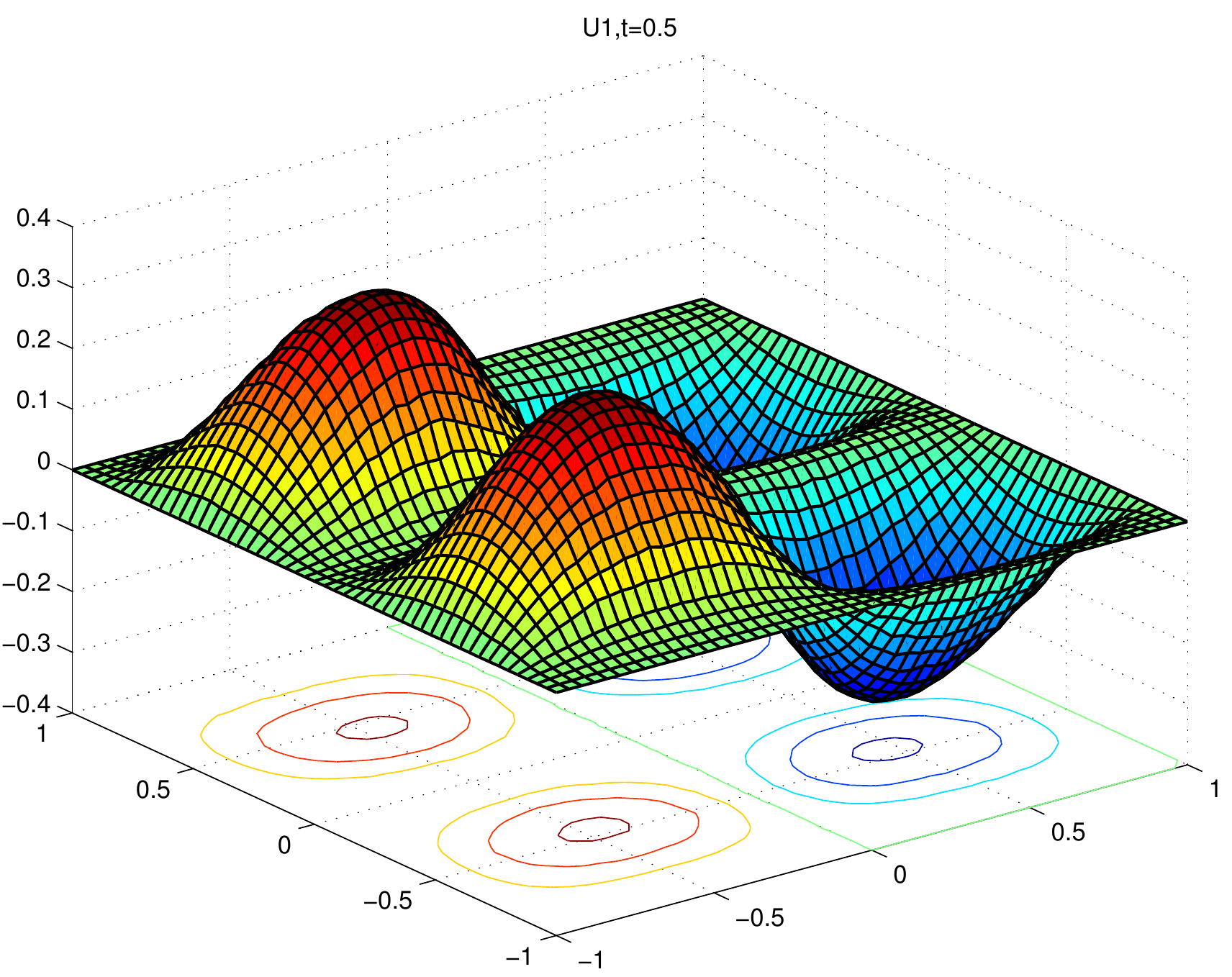}
}\hfill
%\subfigure[$u_2,t=0.5$]{%
     %\includegraphics[width=0.45\linewidth]{figures_finales/Heat2D/M1ref.pdf}
%     \includegraphics[width=0.35\linewidth]{U2T1.pdf}
%}\hfill
\subfigure[$u_1,t=1$]{%
     \includegraphics[width=0.35\linewidth]{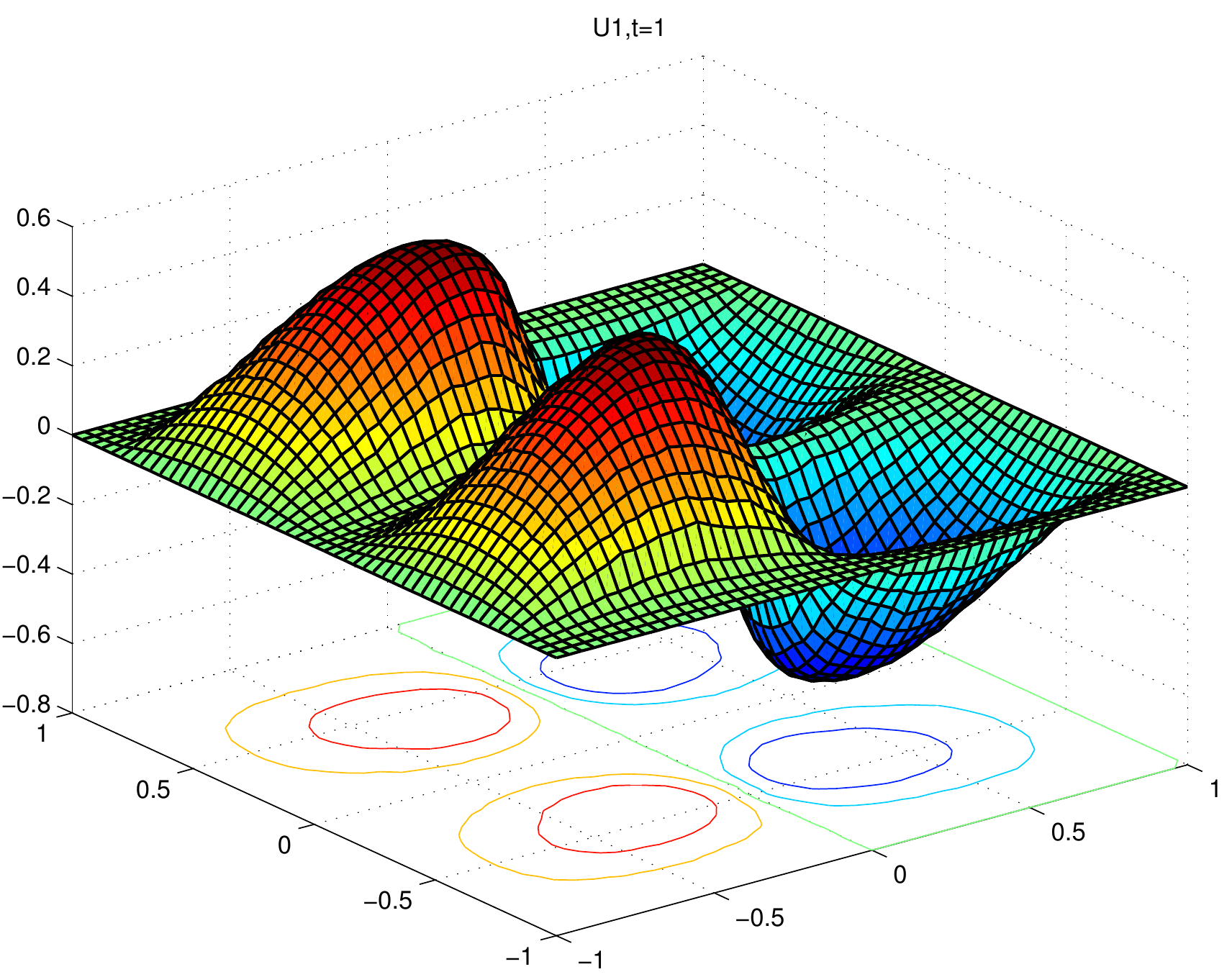}
}\hfill
%\subfigure[$u_2,t=1$]{%
     %\includegraphics[width=0.45\linewidth]{figures_finales/Heat2D/M1.pdf}
%     \includegraphics[width=0.35\linewidth]{U2T2.pdf}
%}\hfill
\subfigure[$u_1,t=1.5$]{%
     \includegraphics[width=0.35\linewidth]{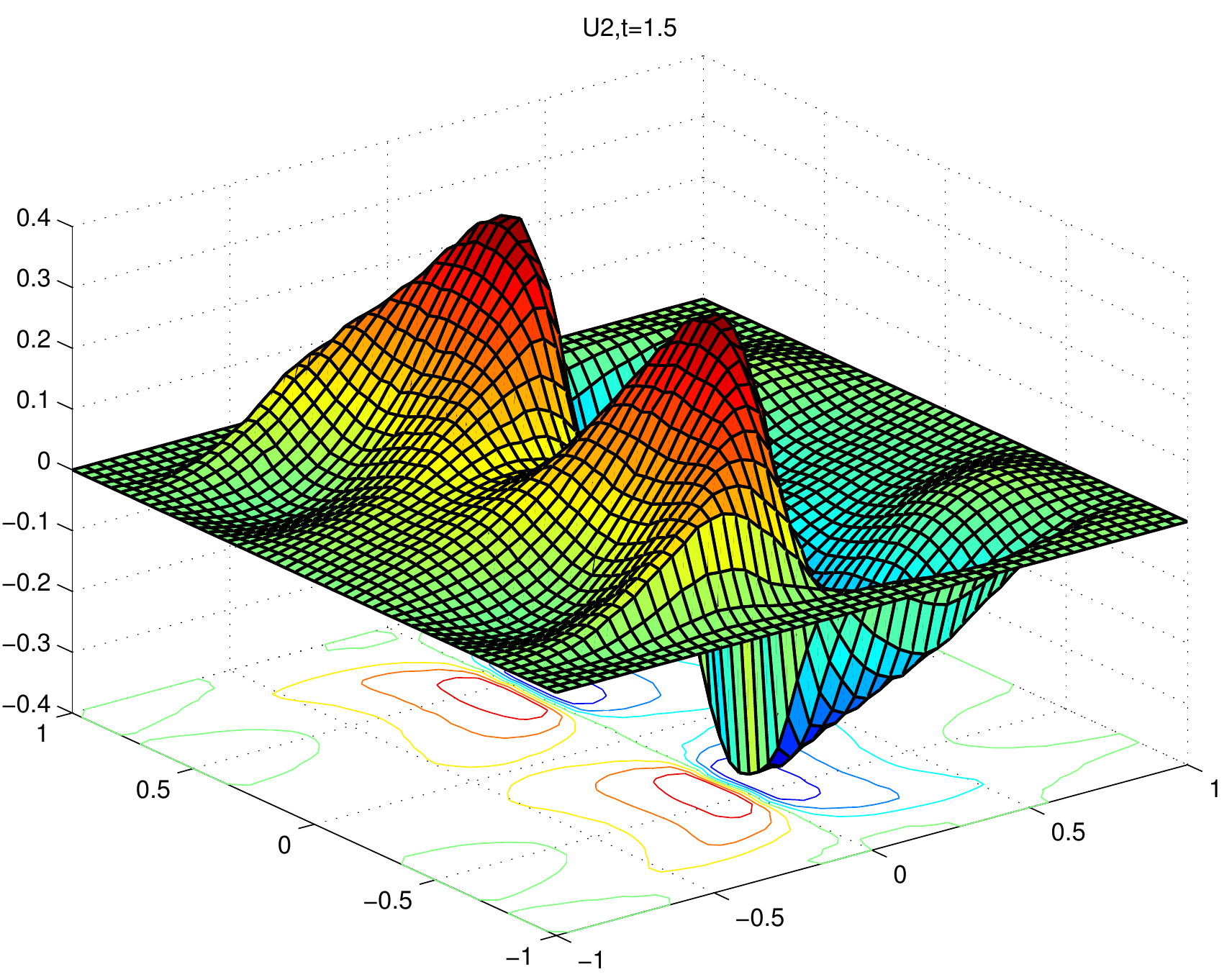}
}\hfill
%\subfigure[$u_2,t=1.5$]{%
     %\includegraphics[width=0.45\linewidth]{figures_finales/Heat2D/M1ref.pdf}
%     \includegraphics[width=0.35\linewidth]{U2T3.pdf}
%}\hfill
\subfigure[$u_1,t=2$]{%
     \includegraphics[width=0.35\linewidth]{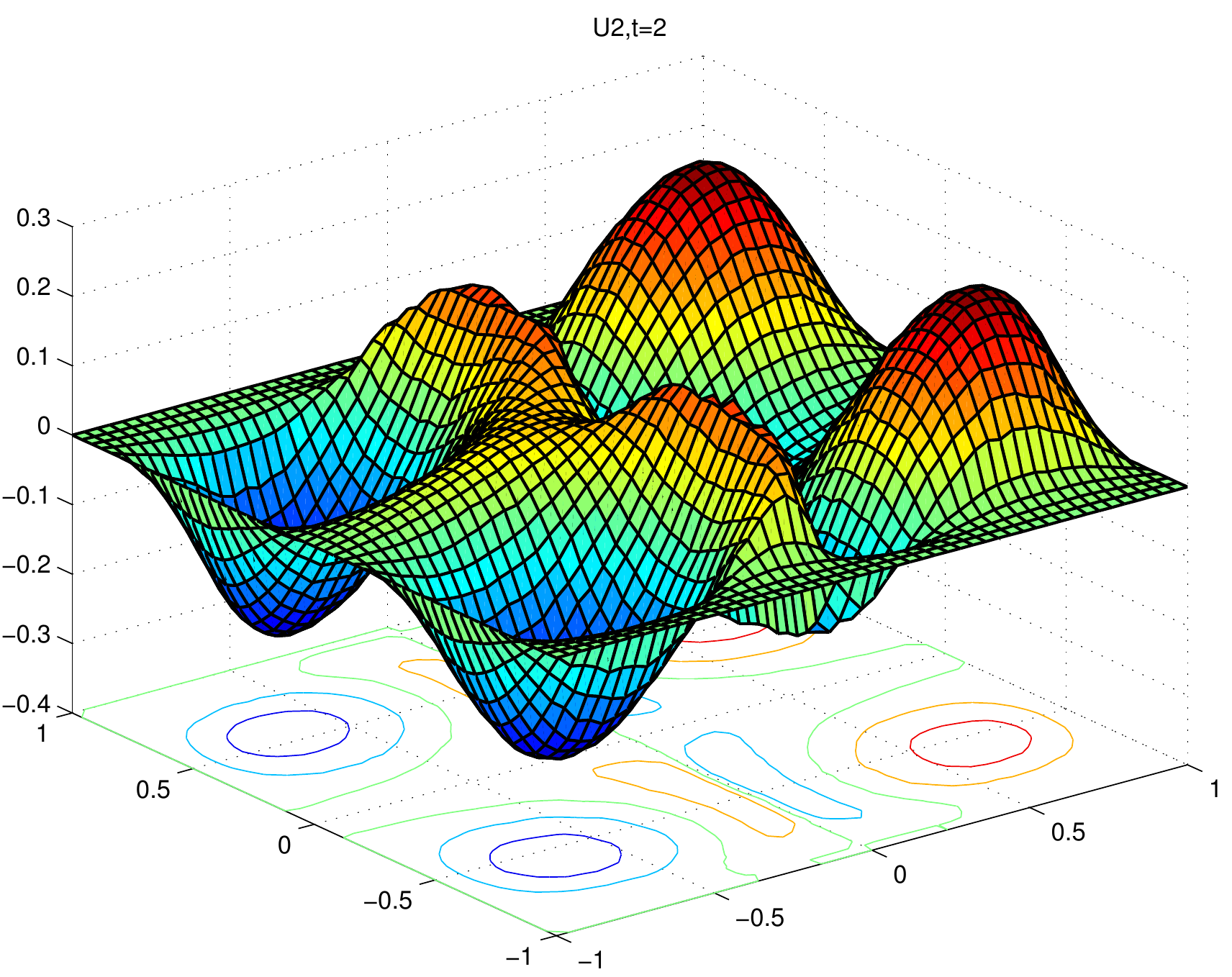}
}\hfill
%\subfigure[$u_2,t=2$]{%
     %\includegraphics[width=0.45\linewidth]{figures_finales/Heat2D/M1ref.pdf}
%     \includegraphics[width=0.35\linewidth]{U2T4.pdf}
%}
\caption{Solution of the 2D Burgers equation at different times - first component}
\label{Burgers2Du1}
\end{figure}

\begin{figure}[h!]
\centering
%\subfigure[$u_1,t=0.5$]{%
     %\includegraphics[width=0.45\linewidth]{figures_finales/Heat2D/M1.pdf}
%     \includegraphics[width=0.35\linewidth]{U1T1.pdf}
%}\hfill
\subfigure[$u_2,t=0.5$]{%
     \includegraphics[width=0.35\linewidth]{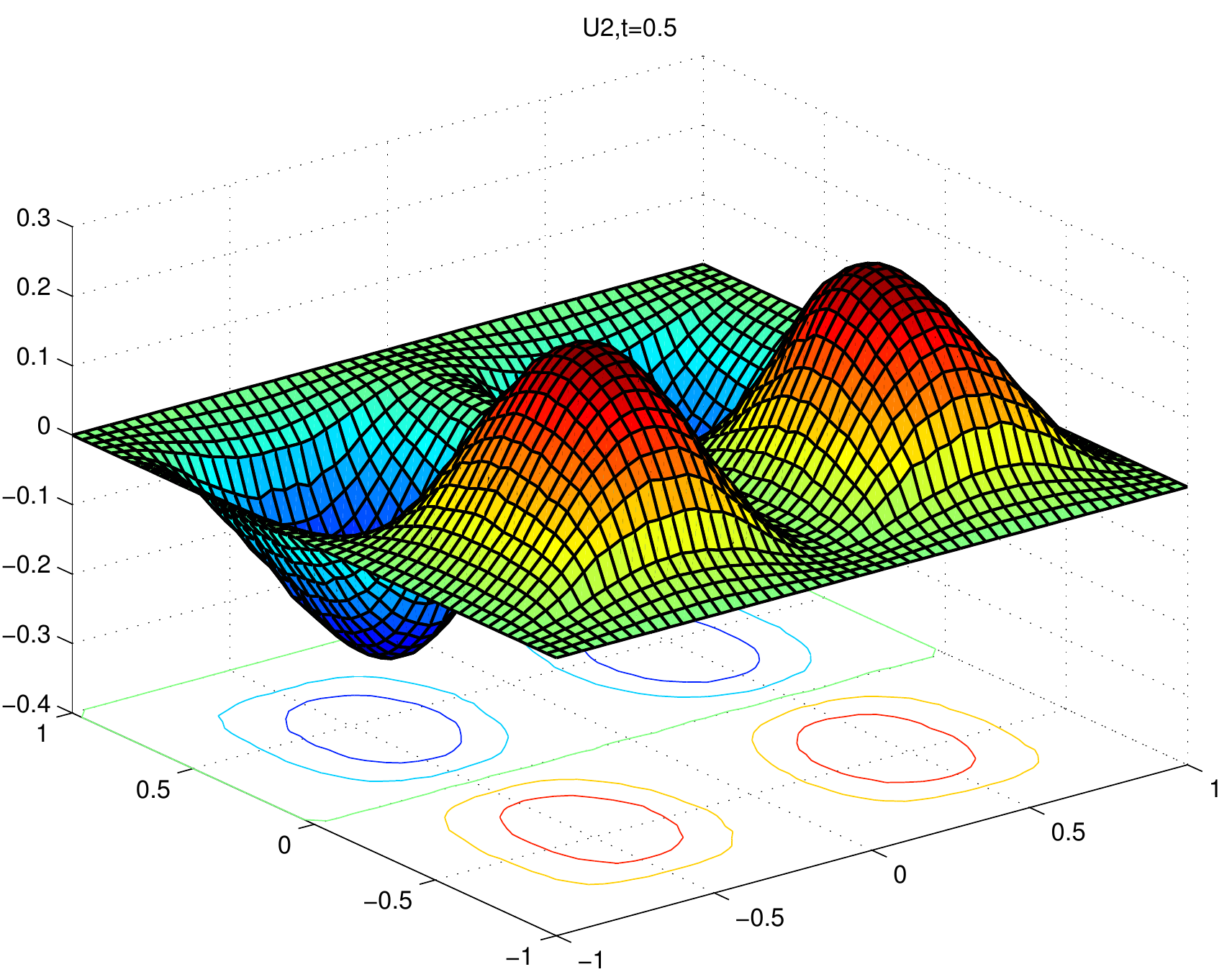}
}\hfill
%\subfigure[$u_1,t=1$]{%
     %\includegraphics[width=0.45\linewidth]{figures_finales/Heat2D/M1.pdf}
%     \includegraphics[width=0.35\linewidth]{U1T2.pdf}
%}\hfill
\subfigure[$u_2,t=1$]{%
     \includegraphics[width=0.35\linewidth]{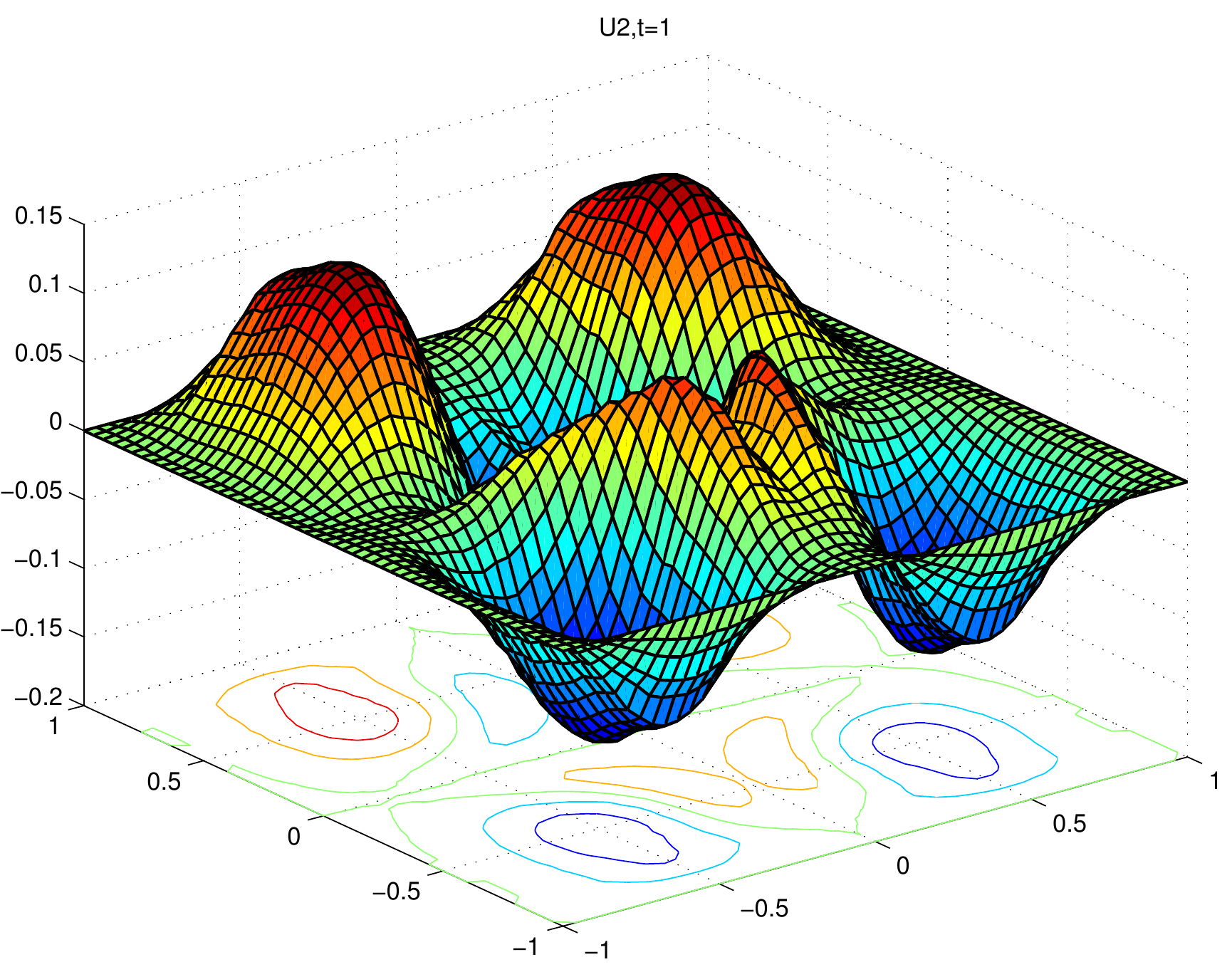}
}\hfill
%\subfigure[$u_1,t=1.5$]{%
     %\includegraphics[width=0.45\linewidth]{figures_finales/Heat2D/M1.pdf}
 %    \includegraphics[width=0.35\linewidth]{U1T3.pdf}
%}\hfill
\subfigure[$u_2,t=1.5$]{%
     \includegraphics[width=0.35\linewidth]{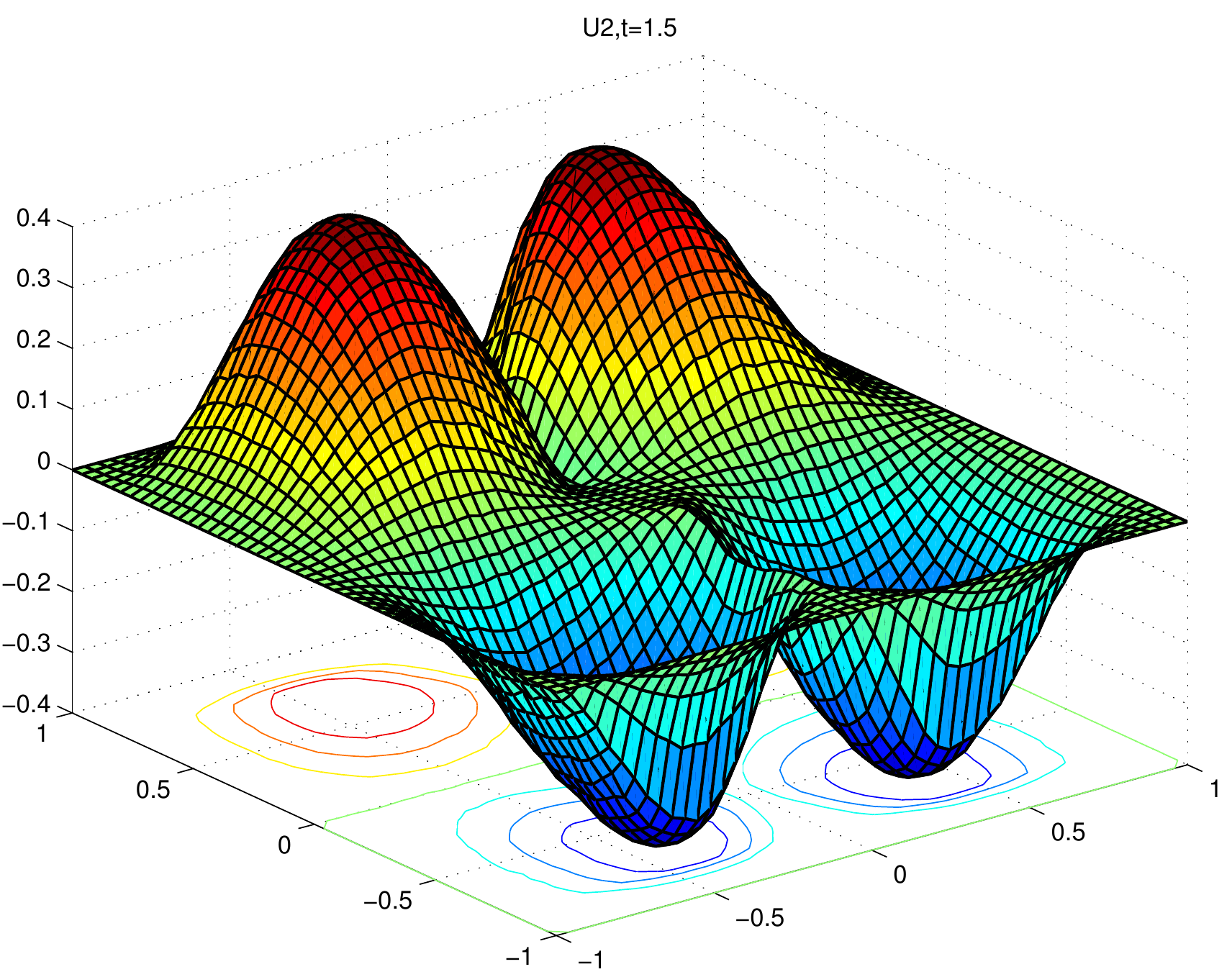}
}\hfill
%\subfigure[$u_1,t=2$]{%
     %\includegraphics[width=0.45\linewidth]{figures_finales/Heat2D/M1ref.pdf}
%     \includegraphics[width=0.35\linewidth]{U1T4.pdf}
%}\hfill
\subfigure[$u_2,t=2$]{%
     \includegraphics[width=0.35\linewidth]{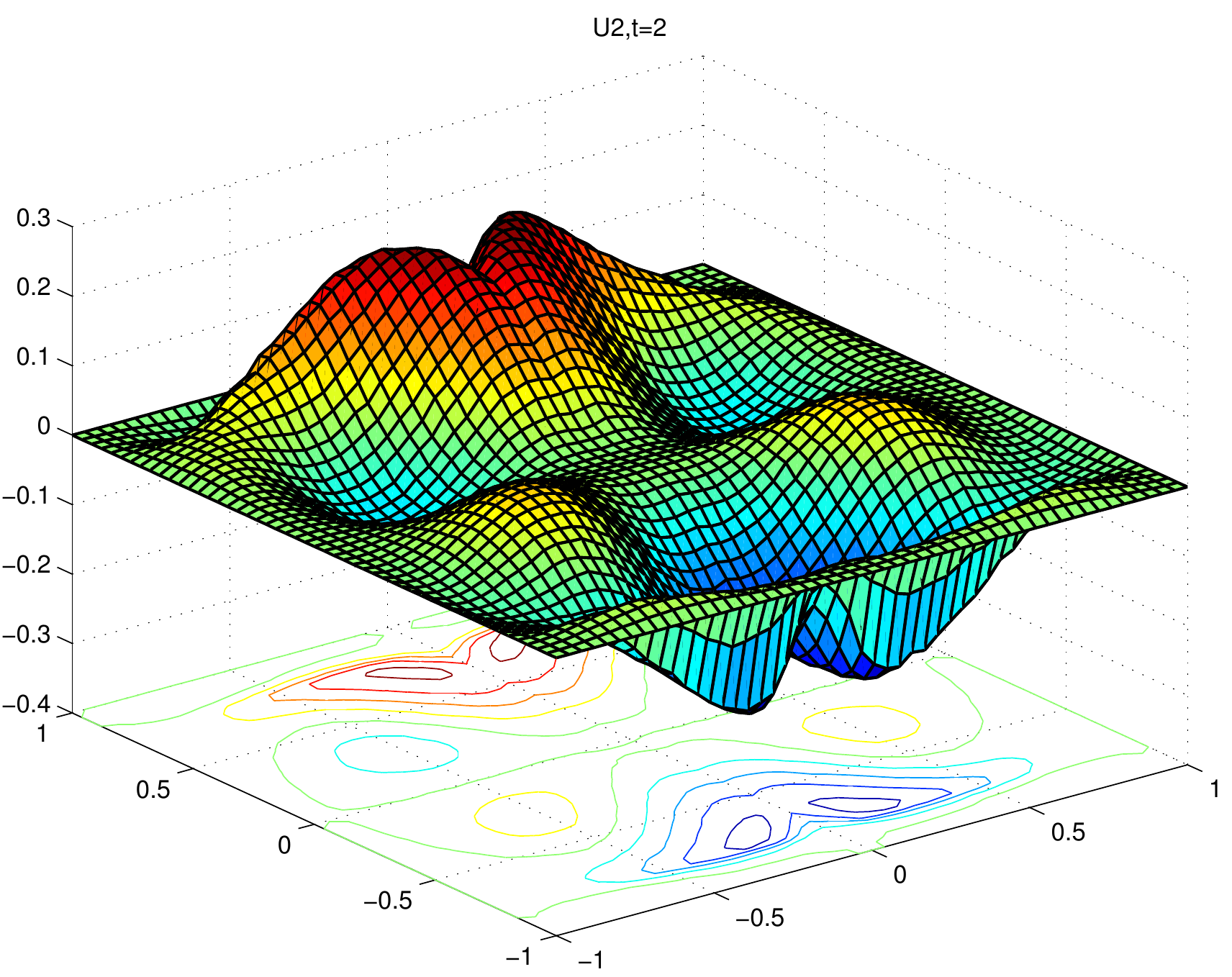}
}
\caption{Solution of the 2D Burgers equation at different times - second component}
\label{Burgers2Du2}
\end{figure}

Both components of the velocity field $u$ are represented in Figures \ref{Burgers2Du1} and \ref{Burgers2Du2} below at different times $t=0.5,1,1.5,2$.

\bibliographystyle{plain}
%\bibliography{bib}

\begin{thebibliography}{1}







\bibitem{crouch}
P. E. Crouch and R. Grossman
{\em Numerical integration of ordinary differential equations on manifolds}, 
J. Nonlinear Sci. 3, 1--33 (1993). 

\bibitem{crouze}
N. Crouseilles, M. Mehrenberger and E. Sonnendr\"ucker, 
{\em Conservative semi-Lagrangian schemes for the Vlasov equation,}
 J. Comput. Phys., 229, 1927--1953 (2010). 

\bibitem{faou}
E. Faou.
\newblock {Analysis of splitting methods for reaction-diffusion problems using
  stochastic calculus.}
\newblock {\em Math. Comput.}, 78(267):1467--1483, 2009.

\bibitem{Falcone}
M. Falcone and R.Ferreti.
{\em Convergence analysis for a class of high-order semi-Lagrangian advection schemes},
SIAM J. Numer. Anal. 35(3), 909 (1998).

\bibitem{Ferretti}
R. Ferretti.
\newblock A technique for high-order treatment of diffusion terms in
  semi-lagrangian schemes, 2000.

\bibitem{gobet}
E. Gobet.
\newblock{Euler schemes and half-space approximation for the simulation of
    diffusion in a domain.}
\newblock{\em ESAIM, Probab. Stat.} 5, 261--297, 2001.

\bibitem{JLL}
B. Jourdain, C. Le~Bris and T.Leli\`evre.
\newblock{On a variance reduction technique for the micro-macro simulations of polymeric fluids.}
\newblock{\em J. Non-Newton. Fluid Mech.}, 122(1-3): 91--106, 2004.


\bibitem{kloeden-platten} {\rm P. Kloeden and E. Platen}, 
{\em Numerical Solution of Stochastic Differential Equations}, Applications of Mathematics (New York) 23, Springer-Verlag, Berlin, 1992.


\bibitem{mann}
R. Mannella.
\newblock{Absorbing boundaries and optimal stopping in a stochastic differential equation.}
\newblock{\em Phys. Lett.,A} 254 (5), 257--262, 1999.

\bibitem{Mil-Tre}
G.N. Milstein and M.V. Tretyakov.
\newblock {\em {Stochastic numerics for mathematical physics.}}
\newblock {Scientific Computation. Berlin: Springer. ixx, 594~p.}, 2004.


\bibitem{Sonnen}
E. Sonnendr\"ucker, J. Roche, P. Bertrand and A. Ghizzo. 
{\em The semi-lagrangian method for the numerical resolution of the Vlasov equation}
J. Comput. Phys. 149, 201--220 (1999)

\bibitem{Staniforth}
A. Staniforth and J. C\^ot\'e, 
\newblock {\em Semi-Lagrangian integration schemes for atmospheric models - A review,}
Mon. Weather Rev. 119 (1991).

\bibitem{Talay96}
{\rm D. Talay}. 
{\em Probabilistic numerical methods for partial differential equations: elements of analysis.} In D. Talay and L. Tubaro (Eds.), Probabilistic Models for Nonlinear Partial Differential Equations, Lecture Notes in Mathematics 1627 (1996) 48--196. 

\end{thebibliography}

\end{document}